\documentclass[11pt,a4paper,reqno]{amsart}
\usepackage{amssymb,amsmath,amsthm}
\usepackage[T1]{fontenc}

\usepackage{color}
\usepackage{enumerate}

\newcommand{\X}{\mathbb{X}} 

\newtheorem{hypo1}{Hypothesis}[section]

\newtheorem{prop1}[hypo1]{Proposition}

\newtheorem{lem1}[hypo1]{Lemma}

\newtheorem{defi1}[hypo1]{Definition}

\newtheorem{rqe1}[hypo1]{Remark}

\newtheorem{coro1}[hypo1]{Corollary}

\newcommand{\Z}{{Z\!\!\!Z}}

\newtheorem{exa}{Application}

\def\B{\mathcal{B}}

\def\PP{\mathbb{P}}
\def\RR{\mathbb{R}}
\def\ZZ{\mathbb{Z}}
\def\CC{\mathbb{C}}
\def\EE{\mathbb{E}}
\def\NN{\mathbb{N}}
\def\SS{\mathbb{S}}

\title{On the recurrence set of planar Markov Random Walks}
       \date{\today }
\author{Lo\"{\i}c Herv\'e}
\address{Universit\'e Europ\'eenne de Bretagne, I.R.M.A.R. (UMR-CNRS 6625), 
Institut National des Sciences Appliqu\'ees
de Rennes, France}
\email{ Loic.Herve@insa-rennes.fr}

\author{Fran\c{c}oise P\`ene}
\address{Universit\'e Europ\'eenne de Bretagne, Universit\'e de Brest,
Laboratoire de Math\'ematiques, UMR CNRS 6205, 29238 Brest cedex, France}
\email{francoise.pene@univ-brest.fr}

\begin{document}
\begin{abstract}
In this paper, we investigate properties of
recurrent planar Markov random walks. 
More precisely, we study the set of recurrence points with the use
of local limit theorems. The Nagaev-Guivarc'h spectral method provides 
several examples for which these local limit theorems are satisfied as soon 
as some (standard or non-standard) central limit theorem and some non-sublattice 
assumption hold.
\end{abstract}

   \maketitle

\noindent{\bf MSC (2010)~:} 60F05

\noindent{\bf Keywords~:} recurrence, Markov chain, spectral method

\section*{introduction} \label{intro1}
Let $\X$ be a measurable space endowed with a $\sigma$-algebra 
${\mathcal X}$. 
Let $(X_n,S_n)_{n\in\NN}$ be a Markov random walk (MRW) with state space 
$\X\times \RR^2$, i.e.~a Markov chain such that
the distribution of $(X_{n+1},S_{n+1}-S_n)$ depends on the past 
only through $X_n$.
Namely: $(X_n,S_n)_{n\in \NN}$ is a Markov chain with transition kernel $P$ satisfying, for any set $A\in {\mathcal X}$ and any Borel subset $S$ of $\RR^2$, the following additive property (in the second component): 
\begin{equation} \label{Add_pro}
\forall (x,s)\in\X\times \RR^2,\ \ P\big((x,s);A\times S\big) = P\big((x,0);A\times (S-s)\big).
\end{equation}
From this definition, the first component $(X_n)_{n\in \NN}$ is 
a Markov chain, called the driving Markov chain of the MRW. 
We suppose that $S_0=0$.  
Given any distribution $\mu$ on $(\X,{\mathcal X})$ (corresponding
to the distribution of $X_0$),  
notation ${\mathbb P}_{(\mu,0)}$ refers to 
the distribution of $(X_n,S_n)_{n\in\NN}$ with initial distribution
$\mu\otimes\delta_0$. 
This notation takes the usual sense when $(X_n,S_n)_{n\in\NN}$ 
is the canonical version defined on $(\X\times \RR^2)^{\NN}$. 
In this work, the last assumption may be assumed without 
loss of generality. The transition kernel of $(X_n)_{n\in\NN}$ is 
denoted by $Q$.

Throughout the paper, we assume that $Q$ admits an invariant probability measure on $\X$, called $\pi$, and that $S_1$ is ${\mathbb P}_{(\pi,0)}$ centered, namely: $S_1$ is 
${\mathbb P}_{(\pi,0)}$-integrable and ${\mathbb E}_{(\pi,0)}[S_1]=0$. Moreover we suppose that there exists a two-dimensional closed subgroup $\SS$ in $\RR^2$ such that we have 
\begin{equation} \label{Sn-supported-by-H}
\forall x\in\X,\  \forall n\in\NN,\ \ \PP_{(x,0)}\big(S_n \in \SS\big) = 1. 
\end{equation}

Let $|\cdot|$ denote the euclidean norm on $\RR^2$. Let us recall that $(S_n)_n$ is said to be recurrent if
$\forall\varepsilon>0,\ {\mathbb P}_\pi\big(|S_n| < \varepsilon\ \text{i.o.}\big) = 1$, 
with the usual notation $[A_n\  \text{i.o.}]:=\cap_{n\ge 1}\cup_{k\ge n} A_k$
("i.o." meaning infinitely often). 
Contrarily to the one-dimensional case,  
the strong law of large numbers (i.e. $S_n/n\rightarrow 0$ a.s.)  
is not sufficient in dimension $2$ to obtain the recurrence property for $(S_n)_n$. 
This is true even in the independent case (which is a special instance of MRW): if $(X_n)_n$ is a sequence of 
$\RR^2$-valued 
independent identically distributed (i.i.d.) centered random variables~(r.v.), then 
$S_n = X_1+\ldots+X_n$ is recurrent if and only if $\sum_n{\mathbb P}(\vert S_n\vert<\varepsilon)=\infty$
for every $\varepsilon>0$.
Hence, in the i.i.d. case, if the distribution of 
$X_1$ is in the domain of attraction of
a stable distribution of index $\alpha$, 
then $\alpha=2$ is required. In other words, in this case, 
a central limit theorem (CLT)   
with a good normalization is needed for $S_n = X_1+\ldots+X_n$ to be recurrent, see \cite[Sect.~3.2]{dur}. 

Recurrence property of $d$-dimensional random walks 
is investigated in many papers. This study is well-known 
for i.i.d.~increments, 
see for instance \cite{dur}. In the dependent case, let us mention 
in particular \cite{berbee,dek,jean-pierre,schmidt-note,schmidt-zwvg,schmidt-ims} 
for random walks with 
stationary increments, \cite{hen1} for MRWs (case $d\geq3$) 
associated with uniformly ergodic 
Markov chains, \cite{gui-local} for MRWs associated with strongly 
ergodic Markov chains,
and \cite{hognas} for additive functionals of Harris recurrent 
Markov chains. 

For general stationary $\RR^2$-valued random walks, the link between 
CLT and recurrence of $(S_n)_n$ 
has been investigated by Conze in \cite{jean-pierre} and
by Schmidt in \cite{schmidt-note} in the situation when the CLT holds 
with the standard normalization in $\sqrt{n}$. 
The methods used in these two works do not extend directly
to other normalizations.

Transience/recurrence properties of MRWs (with $\RR^d$-valued second 
component) have been investigated in \cite{gui-local} 
on the basis of a local limit theorem (LLT) obtained via the standard
Nagaev-Guivarc'h spectral method. This method combined
with the use of the Kochen-Stone
adaptation of the Borel-Cantelli lemma (see (\ref{LLT-1-kochen}) 
and (\ref{LLT-2-kochen}) below) has also been used by Sz\'asz and Varj\'u 
for particular planar stationary walks with standard normalization
in $\sqrt{n}$ in \cite{sza-var} as well as with a non-standard
normalization in $\sqrt{n\log n}$ in \cite{sza-var2}.

Our work uses a similar approach to that of \cite{gui-local,sza-var,sza-var2}, but it goes beyond the question of recurrence. In fact we want to investigate the set of recurrence points ${\mathcal R}_{(\mu,0)}$,
also called {\it recurrence set}, defined by
\[
{\mathcal R}_{(\mu,0)} := \bigg\{\, s\in \SS \, :\, 
\forall\varepsilon>0,\ 
{\mathbb P}_{(\mu,0)}\big(\, |S_n-s| < \varepsilon\ \text{i.o.}\big) = 1\, \bigg\}.  
\]
We simply write ${\mathcal R}_{(x,0)}$ when $\mu$ is the Dirac distribution $\delta_x$ at $x\in\X$.
We describe situations in which we prove that ${\mathcal R}_{(\mu,0)} = \SS$  
for every initial distribution $\mu$.

The recurrence set is well-known in the i.i.d.~case 
(e.g.~see \cite[Sect.~3.2]{dur}), 
and it has been fully investigated in \cite{als} for 
one-dimensional MRW (i.e.~$S_n$ is real-valued). However, to the best 
of our knowledge, the recurrence set has not been 
investigated for planar MRWs.

In dimension 2, whereas some recurrence results are only based on the CLT,
the study of the recurrence set requires some assumption 
ensuring (roughly speaking) 
that $\SS$ is the smallest lattice in $\RR^2$ 
satisfying (\ref{Sn-supported-by-H}). 
Note that such a lattice-type assumption is also the additional 
condition to pass from CLT to LLT. 
Therefore, it is not surprising that local limit theorems will play here an important role in
the study of the set ${\mathcal R}_{(\mu,0)}$ for planar MRWs. 
The LLTs involved in this work are obtained by using the weak Nagaev-Guivarc'h
spectral method developed in \cite{loicsoaz}. 

In Section \ref{presentation}, we state our main results and give applications.
First we state two key theorems (Theorems~I-II) giving 
${\mathcal R}_{(\mu,0)}=\mathbb S$ under conditions related to LLTs.
Second we state two theorems (Theorems~III-IV) 
giving these local limit conditions, under general assumptions, 
as soon as $(S_n)_n$ satisfies a (standard or non-standard) 
central limit theorem with suitable normalization
as well as a non-sublattice condition in $\mathbb S$.
We illustrate our general results with classes of models.

In Section \ref{simp-rks-S}, we make some simple remarks on the subgroup 
$\mathbb S$ appearing in (\ref{Sn-supported-by-H}), on the non-sublattice condition in $\mathbb S$, and on the special case when $S_n$ is an additive functional of $(X_n)_{n\in \NN}$.  

In Section \ref{sec:basic-fact}, we prove Theorems~I-II by using classical arguments 
derived from the i.i.d.~case \cite{dur} and the Kochen and Stone adaptation 
of the Borel-Cantelli lemma \cite{KochenStone}. 
In Sections \ref{sec-meth-spec} and \ref{sec-proof-deux-prop}, 
we prove Theorems~III-IV thanks to the weak Nagaev-Guivarc'h
spectral method and Fourier techniques. In Section \ref{sect-ex-1}, we detail our applications (that have been
shortly introduced in Section~\ref{presentation}). These applications
are obtained thanks to recent works \cite{BDG,djl,jvl,loicsoaz}. Some complements concerning spectral method are given in appendices. 

\noindent{\bf Acknowledgments.--}  {\it We wish to thank the referee for his useful advices which allowed us to improve the presentation of our article.} 

\section{Main results and applications} \label{presentation}
Let $(X_n,S_n)_{n\in \NN}$ be a Markov random walk with state space $\X\times\RR^2$, 
let $\SS$ be a two-dimensional closed subgroup of $\RR^2$ 
satisfying (\ref{Sn-supported-by-H}). The Haar measure on $\SS$ is called $m_\SS$. 
We denote by $B(s,\varepsilon)$ the open ball in $\RR^2$  
centered at $s$ with radius $\varepsilon$. We denote by $\B(\RR^2)$ 
the Borel $\sigma$-algebra of  $\RR^2$. 

We start by stating two key results giving two different approaches to prove 
${\mathcal R}_{(\mu,0)}=\mathbb S$.
The first theorem borrows classical arguments 
derived from the i.i.d.~case \cite{dur}. Recall that $(X_n)_{n\in\NN}$ is said to be 
Harris recurrent if, for any set $B\in{\mathcal X}$ such that $\pi(B)>0$, 
for every $x\in\X$, we have 
${\mathbb P}_{(x,0)}\big(X_k\in B\ \text{i.o.}\big) = 1$. 
\smallskip

\noindent{\bf Theorem~I.} {\it 
Assume that there exist $\varepsilon_\SS>0$ and a 
sequence $(a_n)_{n\ge 1}$ of positive real numbers satisfying 
$\sum_{n\ge 1}a_n=\infty$ such that, for every 
$(s,\varepsilon)\in \SS\times(0;\varepsilon_\SS)$, 
for every  bounded measurable function $f : \X\rightarrow [0,+\infty)$, the following local limit property holds with $B:=B(s,\varepsilon)$: 
\begin{equation} \label{LLT-f-bornee}
\EE_{(\pi,0)}\big[f(X_n)\, {\mathbf 1}_B(S_n)\big] \sim a_n\, \pi(f)\, m_\SS(B)\quad \text{when } n\rightarrow+\infty. 
\tag{\text{LL0}}
\end{equation} 
Then the following assertions hold true:
\begin{itemize}
\item[\bf (a)] ${\mathcal R}_{(\pi,0)} = \SS$; 
\item[\bf (b)] if in addition $(X_n)_{n\in\NN}$ is Harris recurrent, 
then ${\mathcal R}_{(\mu,0)}=\SS$ for every  
probability measure $\mu$ on $(\X,{\mathcal X})$. 
\end{itemize}
} 
\smallskip

Theorem I is direct and quite natural but it requires a local limit estimate 
for every bounded measurable functions. Moreover it needs Harris-recurrence
hypothesis to obtain the non-stationary result {\bf (b)}.
These two assumptions appear to be not satisfied on some classes of models.
For this reason, we give another theorem based on another approach
but still related to local limit theorems.
Theorem~II below involves the Kochen and Stone adaptation 
of the Borel-Cantelli lemma.  
\smallskip

\noindent{\bf Theorem~II.} {\it 
Assume that there exist $\varepsilon_\SS>0$ and $d>0$ 
such that, for every $\varepsilon\in(0;\varepsilon_\SS)$, 
for every $(x,s)\in \X\times \SS$, we have 
\begin{equation} \label{cond-cor-resume-1}
 \sum_{n\ge 1}{\mathbb P}_{(x,0)}
\big(\vert S_n-s\vert<\varepsilon\big) = \infty, 
 \ \ \quad \quad \quad \quad \quad\quad \quad \quad \quad \quad \quad
 \quad \quad \quad \quad 
\tag{\text{KSa}}
\end{equation}
\begin{equation} \label{cond-cor-resume-2}  
\liminf_{N\rightarrow +\infty} \frac{ \sum_{n,m=1}^N{\mathbb P}_{(x,0)}\big( \vert S_n-s\vert<\varepsilon,\, \vert S_{n+m}-s\vert
<\varepsilon\big)}
{\left(\sum_{n=1}^N{\mathbb P}_{(x,0)}( \vert S_n-s\vert<\varepsilon)  \right)^2} \leq d. \quad \quad \quad \quad \quad 
\tag{\text{KSb}}
\end{equation}
Then ${\mathcal R}_{(\mu,0)} = \SS$ for every 
probability measure $\mu$ on $(\X,{\mathcal X})$. }
\smallskip

It is not difficult to prove (see Lemma~\ref{sim-sum}) that the 
Kochen-Stone Conditions~(\ref{cond-cor-resume-1})-(\ref{cond-cor-resume-2}) 
are implied by the two following local limit properties
(the second being a bivariate local limit theorem): 
\begin{equation} \label{LLT-1-kochen}
{\mathbb P}_{(x,0)} \big(S_n\in B\big) \sim D\, a_n m_\SS(B), \quad \quad \quad \quad \quad 
\tag{\text{LLa}}
\end{equation}
\begin{equation} \label{LLT-2-kochen} 
{\mathbb P}_{(x,0)}\big((S_n,S_{n+m})\in B^2\big) \sim D^2\, a_na_m\,  m_\SS(B)^2, \quad \quad \quad \quad \quad 
\tag{\text{LLb}}
\end{equation}
where $B:=B(s,\varepsilon)$, $D\in(0,+\infty)$ and where $(a_n)_{n\ge 1}$ 
is a sequence of positive numbers satisfying $\sum_{n\ge 1}a_n=\infty$. 

Now let us present, and illustrate by classes of models, our general operator-type strategy 
providing (\ref{LLT-f-bornee}) and 
(\ref{LLT-1-kochen})-(\ref{LLT-2-kochen}), and so conclusions of Theorems~I and II. 
To that effect, 
we consider the Fourier operators $Q(t)\, $ ($t\in\RR^2$), associated with the MRW, 
acting (in a first step) on the space of 
bounded measurable functions $f : \X\rightarrow{\mathbb C}$, as follows: 
\begin{equation} \label{Fourier}
\forall t\in\RR^2,\ \forall x\in \X,\  \ 
\big(Q(t)f\big)(x) := \EE_{(x,0)}\big[e^{i \langle t , S_1 \rangle} f(X_1)\big]. 
\end{equation}

If $({\mathcal Y},\|\cdot\|_{{\mathcal Y}})$ is a Banach space, 
${\mathcal L}({\mathcal Y})$ denotes the space of linear continuous 
endomorphisms of $\mathcal Y$. The associated operator norm is also 
denoted by $\Vert\cdot\Vert_{\mathcal Y}$. 

Let $({\mathcal B},\|\cdot\|_{{\mathcal B}})$ and $(\widehat{\mathcal B},\|\cdot\|_{\widehat{\mathcal B}})$ be two complex Banach spaces composed 
of $\pi$-integrable $\CC$-valued functions on $\X$ (or of classes 
modulo $\pi$ of such functions). 
We assume that ${\mathcal B} \subset \widehat{\mathcal B}$ and that 
\begin{equation} \label{cont-inclusions}
\exists c\in(0,+\infty),\ \forall f\in\mathcal B,\ \|f\|_{\widehat{\mathcal B}} \leq c\|f\|_{\mathcal B} 
\quad \text{and}\quad 
\exists d\in(0,+\infty),\ \forall f\in\widehat{\mathcal B}, \ \int|f|d\pi \leq d\|f\|_{\widehat{\mathcal B}}.
\end{equation}
For ${\mathcal Y} = {\mathcal B}$ or 
$\widehat{\mathcal B}$ and for any probability measure $\nu$ on $\X$, we write $\nu\in{\mathcal Y}'$ when the linear map $f\mapsto \nu(f) := \int_{\X} f\, d\nu$ is well-defined and continuous on $\mathcal Y$. 
The conditions in (\ref{cont-inclusions}) imply that 
$\pi\in {\mathcal B}'\cap(\widehat{\mathcal B})'$. 
We denote by $\Vert \cdot\Vert_{{\mathcal B},\, \widehat{\mathcal B}}$ the operator 
norm of bounded linear operators from ${\mathcal B}$ to $\widehat{\mathcal B}$. 

The essential spectral radius of any $T\in{\mathcal L}({\mathcal B})$ is defined by $r_{ess}(T) := \lim_n(\inf \|T^n-K\|_{\mathcal B})^{1/n}$ where the infimum is taken over the ideal of compact endomorphisms $K$ on  $\mathcal B$, see e.g.~\cite{wu}. The following assumptions will be used to check hypotheses of Theorems~I-II. 
\medskip

\noindent{\bf Operator-type assumptions.} {\it 
Function ${\bf 1}_\X$ (or its class) is in $\mathcal B$; for every $t\in {\mathbb R}^2$ we have 
$Q(t)\in{\mathcal L}({\mathcal B})\cap{\mathcal L}(\widehat{\mathcal B})$.
Moreover: 
\begin{itemize} \label{HYPO}
  \item[\bf (A1)] $Q$ is strongly ergodic on ${\mathcal B}$, 
that is: $\lim_n Q^n = \pi$ in ${\mathcal L}({\mathcal B})$,     
  \item[\bf (A2)] $\forall t\in{\mathbb R}^2,\quad \lim_{h\rightarrow 0}
     \Vert Q(t+h)-Q(t)\Vert_{{\mathcal B},\, \widehat{\mathcal B}}=0$,  
  \item[{\bf (A3)}] for every compact $K$ in ${\mathbb R}^2$, there exist $\kappa\in(0;1)$, 
$C\in(0;+\infty)$ such that, for every $t\in K$, 
the essential spectral radius of $Q(t)$ on $\mathcal B$ satisfies $r_{ess}(Q(t))\le\kappa$, and  
\begin{equation} \label{ineg-DF} 
\forall n\ge 1,\ \forall f\in\mathcal B,\quad \Vert Q(t)^nf\Vert_{\mathcal B}\le C\kappa^n\Vert f\Vert_{\mathcal B}
      +C\Vert f\Vert_{\widehat{\mathcal B}},
\end{equation}
  \item[{\bf (A4)}] for every $\lambda\in\mathbb C$ such that $\vert\lambda\vert\ge 1$
and for every nonzero element $f\in\mathcal B$, the implication
$$\left[\exists n_0,\ \forall n\ge n_0,\ |\lambda|^n|f|\le Q^n|f|\right]\ 
\Rightarrow\ \left[\vert\lambda\vert=1\ \mbox{and}\ \vert f\vert\le\pi(|f|)\right]$$
holds true. 
\end{itemize}
}
These operator-type assumptions put all together correspond 
to Conditions~$(\widehat{K})$ and $(P)$ of \cite[p.~434-436]{loicsoaz}. They are the key assumptions to use the weak Nagaev-Guivarc'h method. Further comments on these hypotheses are presented at the end of Section~\ref{sec-meth-spec}. 

\medskip

\noindent{\bf Probabilistic-type assumptions:}  
\begin{itemize} {\it 
  \item[{\bf (A5)}] There exist a 
function $L:(0;+\infty)\rightarrow(0;+\infty)$ assumed to be slowly varying at $\infty$
(i.e.~$\forall k>0$, $\lim_{x\rightarrow +\infty} L(kx)/L(x)=1$) 
and a sequence $(A_n)_n$ of positive real numbers satisfying $A_n^2\sim nL(A_n)$ and $\sum_{n\geq1} A_n^{-2}=\infty$ such that, under ${\mathbb P}_{(\pi,0)}$, $(S_n/A_n)_n$ converges in distribution to a non-degenerate Gaussian law  ${\mathcal N}(0,\Gamma)$, \\

	\item[{\bf (A6)}] $(S_n)_n$ is {\bf non-sublattice in $\SS$}; namely 
there exists no $(\SS_0,\chi,(\beta_t)_{t\in\SS_0^*})$ 
with $\SS_0$ a closed proper subgroup of $\, \SS$, 
$\, \chi : {\mathbb X}\rightarrow{\mathbb R}^2$ a bounded 
measurable function, $(\beta_t)_{t\in \SS_0^*}$ a family 
of real numbers indexed by the dual group $\SS_0^*$ of $\SS_0$, 
satisfying the following property for $\pi$-almost every $x\in\X$: 
\begin{equation} \label{fct-beta} 
\forall t\in\SS_0^*,\ \forall n\geq1,\quad \big\langle t, S_n + \chi(X_n) - \chi(x)\big\rangle \in 
n \beta_t + 2\pi\Z\ \ {\mathbb P}_{(x,0)}-a.s..
\end{equation}
}
\end{itemize}

\medskip

\noindent{\bf Theorem III.} {\it Assume that assumptions (A1) to (A6) hold true
with $\mathcal B$ containing all the nonnegative bounded measurable
functions.
Then (\ref{LLT-f-bornee}) is fulfilled, and so
conclusions of Theorem~I apply. 
}

\medskip

Let us observe that the assumption on $\mathcal B$ in Theorem III is not fulfilled 
if $\mathcal B$ is defined as some space of regular functions. In this case, 
the next statement is relevant, and it is worth noticing that it does not require 
the Harris-recurrence hypothesis.

In Theorem IV below, we suppose that $\widehat{\mathcal B}$ is composed of 
functions (not of classes), so that the Dirac distribution at any $x\in\X$, called $\delta_x$, 
is defined on $\widehat{\mathcal B}$ (i.e.~$\delta_x(f):=f(x)$). Note that 
$\delta_x\in(\widehat{\mathcal B})'$ means that there exists $c_x>0$ such that: 
$\forall f\in\widehat{\mathcal B},\ |f(x)| \leq c_x\, \Vert f\Vert_{\widehat{\mathcal B}}$.

\medskip

\noindent{\bf Theorem IV.} {\it 
Assume that assumptions (A1) to (A6) hold true. Suppose that, for every $x\in\X$,
$\delta_x\in(\widehat{\mathcal B})'$.
Then (\ref{LLT-1-kochen}) and (\ref{LLT-2-kochen}) are fulfilled, thus the
conclusion of Theorem~II holds.}

\medskip

Our method enables the study of recurrence set for every
model satisfying our general assumptions {\bf (A1)}-{\bf (A4)}. The development of the use of the 
Nagaev-Guivarc'h method offers large perspectives of applications.
To fix ideas, we give now particular applications of 
Theorems III and IV to classes of models. 
These applications are proved and detailed in Section \ref{sect-ex-1}. 

The following result will be derived from Theorem~III.
\begin{exa}[$\rho$-mixing driving chain] \label{aplli-rho-mix}
Assume that the driving Markov chain $(X_n)_n$ is $\rho$-mixing, that
$S_1$ is centered and square integrable under ${\mathbb P}_{(\pi,0)}$, that $(S_n)_n$ is non-sublattice in $\SS$ and that the limit covariance matrix $\Gamma$
of $(S_n/\sqrt{n})_n$ is positive definite.
Then ${\mathcal R}_{(\pi,0)}=\SS$. If in addition $(X_n)_{n\in\NN}$ is Harris recurrent, 
then ${\mathcal R}_{(\mu,0)}=\SS$ for any initial distribution $\mu$.
\end{exa}

The three following results will be derived from Theorem~IV. They concern the special case when $S_n$ is defined as a   univariate or bivariate additive functional (AF) of $(X_n)_n$. 

\begin{exa}[AF of $V$-geometrically ergodic Markov chain] \label{appli-v-geo}
Let us assume that $(X_n)_n$ is $V$-geometrically ergodic for some
$V: \X\rightarrow [1,+\infty)$ and that 
\begin{itemize}
\item[(a)] either $S_n=\sum_{k=1}^n\xi(X_k)$ for some $\pi$-centered 
$\xi:\X\rightarrow{\mathbb R}^2$ such that 
$\vert \xi\vert^2/V$ is bounded, \\
\item[(b)] or $S_n=\sum_{k=1}^n\xi(X_{k-1},X_k)$ for some 
$\xi:\X\times\X\rightarrow{\mathbb R}^2$ such that $\xi(X_0,X_1)$ is 
${\mathbb P}_\pi$-centered and 
$\sup_{x,y} \vert\xi(x,y)\vert^{2+\varepsilon}/(V(x)+V(y))$ is finite for some $\varepsilon>0$. 
\end{itemize}
If $(S_n)_n$ is non-sublattice in $\SS$  and the limit covariance matrix $\Gamma$
of $(S_n/\sqrt{n})_n$ is positive definite, then, for every initial distribution 
$\mu$ of $(X_n)_{n\in\NN}$, the recurrence set 
${\mathcal R}_{(\mu,0)}$ of $(S_n)_n$ satisfies 
${\mathcal R}_{(\mu,0)}=\SS$.
\end{exa}
\begin{exa}[AF of Lipschitz iterative models] \label{appli-IFS}
Let $b\ge 0$.
Suppose that $(X_n)$ is a random walk on $\X=\mathbb R^D$ 
given by $X_n = F(X_{n-1},\vartheta_n),\ \ n\geq 1$,
with $F : \X\times V\rightarrow \X$ a measurable function and
with $(\vartheta_n)_{n\geq 1}$ a sequence of i.i.d. random variables 
independent of $X_0$.
Assume that, almost surely, $F(\cdot,\vartheta_1)$ has
Lipschitz constant strictly less than 1, that $ {\mathbb E}[d\big(F(0,\vartheta_1),
0\big)^{2(b+1)}]<\infty$. 
Take
$S_n=\sum_{k=1}^n\xi(X_k)$ for some $\pi$-centered 
$\xi:\X\rightarrow{\mathbb R}^2$. Finally assume 
that $(S_n)_n$ is non-sublattice in $\SS$ and that
\begin{equation} \label{xi-iterative-bis} 
\exists C\geq 0,\ \ \forall (x,y)\in \X^2,\ \
 \big|\xi(x)-\xi(y)\big| \leq C\, d(x,y)\,
 \big[1+d(x,x_0)+d(y,x_0)\big]^b.
\end{equation} 
Then ${\mathcal R}_{(\mu,0)}=\SS$ 
for every initial distribution $\mu$ on $\X$.
\end{exa}

The previous applications involve the standard CLT. In fact it is not so easy to find examples of MRWs, even in case of additive functionals, for which $S_n$ satisfies a non-standard CLT. 
Such instances can be found in \cite{BDG}, see also \cite{mir}. The following application, based on \cite{BDG}, shows how our results apply to affine recursions which are special instances of Lipschitz iterative models. 

\begin{exa}[Affine recursion with non-standard CLT] \label{appli-non-stand}
Suppose that $(X_n)_n$ is a Markov chain on ${\mathbb R}^2$ given by
$X_n=A_nX_{n-1}+B_n,$
where $(B_n,A_n)_n$ is a sequence of i.i.d. 
$\mathbb R^2\rtimes Sim({\mathbb R}^2)$-valued random variables 
($Sim({\mathbb R}^2)$ being the 
similarity group of ${\mathbb R}^2$) independent of $X_0$. 
Assume 
$${\mathbb E}[|A_1|^2]= 1,\quad {\mathbb E}[\vert A_1|^2\log|A_1|]<\infty,\quad {\mathbb E}[|B_1|^2]<\infty.$$
Under some additional conditions (to be specified in Subsection \ref{subsec-aff-rec}) on the 
support of the distribution of $(B_1,A_1)$ and of the invariant probability measure, it is
proved in \cite{BDG} that
there exist $m_0\in{\mathbb R}^2$ and a gaussian random variable $Z$
such that, for every $x\in\mathbb R^2$, $(S_n/\sqrt{n\log(n)})_n$
converges in distribution (under $\mathbb P_x$) to $Z$, with
$S_n:=\sum_{k=1}^n(X_k-m_0)$.

Then, if $(S_n)_n$ is non-sublattice in ${\mathbb R}^2$, 
we have ${\mathcal R}_{(\mu,0)}={\mathbb R}^2$ for every initial measure $\mu$.
\end{exa}

%
%
%
%
\section{Preliminary remarks on Hypotheses~(\ref{Sn-supported-by-H}) and {\bf (A6)}} \label{simp-rks-S}
Given any subgroup $H$ of $\RR^2$, its dual subgroup is defined as 
\begin{equation}\label{H*}
H^*:=\{s\in{\mathbb R}^2\ :\ \forall t\in H,\ \langle s,t\rangle\in 2\pi{\mathbb Z}\}. 
\end{equation}
$H^*$ is a subgroup of $\RR^2$, and the dual subgroup of $H^*$ (i.e.~the bidual of $H$) 
coincides with $H$. These properties are classical, anyway they can 
easily be proved in Cases {\it (H1) (H2) (H3)} below. 

\noindent{\scriptsize $\bullet$} {\it Remarks on (\ref{Sn-supported-by-H}).} 
Theorems I and II are valid for $d$-dimensional MRW. However, in practice, the condition $\sum_{n\ge 1}a_n=\infty$ involved in (\ref{LLT-f-bornee}) and (\ref{LLT-1-kochen})-(\ref{LLT-2-kochen}) is only fulfilled in dimension $d=1$ or~2. 
The one-dimensional cases $\SS=\RR\, \vec u$ and $\SS=\ZZ\, \vec u\, $ ($\vec u\in\RR^2$) are not investigated here since 
the recurrence set of $S_n$ can be deduced from \cite{als} (thanks to the strong law of large numbers).

Consequently, throughout the paper, the subgroups $\SS$ of interest in (\ref{Sn-supported-by-H}) are the two-dimensional closed subgroups of $\RR^2$, which correspond to the three following cases: 
\begin{enumerate} {\it 
\item[(H1)] $\, \SS=\RR^2$. We have $\SS^* = \{0\}$. We set $\varepsilon_\SS := 1$.  \\
\item[(H2)] There exists $(b,\vec u,\vec v)\in(0,+\infty)\times\RR^2\times\RR^2$ 
such that: 
$\SS = b\, \ZZ\, \vec u\, \oplus\, \RR\, \vec v$. We suppose, without loss of 
generality, that $(\vec u,\vec v)$ is an orthonormal basis of $\RR^2$. 
We set $\varepsilon_\SS := b$. Note that, for every $(s,\varepsilon)\in \SS\times(0;\varepsilon_\SS)$, 
we have $B(s,\varepsilon)\cap \SS=\{s+w\vec v\ :\ w\in\RR,\, \vert w\vert<\varepsilon\}$, and that $\SS^* = a\, \ZZ\, \vec u$ with $a=2\pi/b$.  \\
\item[(H3)] There exists some real-valued invertible 
$2\times2-$matrix $B$ such that: $\SS = B\, \ZZ^2$. We set $\varepsilon_\SS := 
\min\{\vert s\vert;\ s\in \SS\setminus \{0\}\}$. Note that, for every 
$(s,\varepsilon)\in \SS\times(0;\varepsilon_\SS)$, $B(s,\varepsilon)\cap \SS=\{s\}$, and that $\SS^* = A\, \ZZ^2$ with $A = 2\pi\, (B^*)^{-1}$, where $B^*$ is the transpose matrix of $B$. }
\end{enumerate}

\noindent{\scriptsize $\bullet$} {\it Remarks on non-sublattice condition {\bf (A6)}.} 
Since the dual subgroup of $\SS_0^*$ is $\SS_0$, one can easily check that, if $(s,s')\in\RR^2$ is such that $\langle t,s\rangle\in n \beta_t + 2\pi\Z$ and $\langle t,s'\rangle\in n \beta_t + 2\pi\Z$ for every $t\in\SS_0^*$, 
then $s-s'\in\SS_0$, namely: $s$ and $s'$ belong to the same class modulo $\SS_0$. Therefore the non-sublattice assumption {\bf (A6)} 
is equivalent to the nonexistence of $(\SS_0,(b_n)_n,\chi)$ with $\SS_0$ 
a proper subgroup  of $\SS$, $(b_n)_n$ a sequence of vectors in $\RR^2$, $\chi:{\mathbb X}\rightarrow\mathbb R^2$ 
a bounded measurable function  such that: 
\begin{equation} \label{gamma-t}
\forall t\in\SS_0^*,\ \exists \beta_t\in\RR,\ \forall n\geq1,\quad \langle t,b_n\rangle\in n 
\beta_t + 2\pi\Z
\end{equation} 
and such that, for $\pi$-almost every $x\in\mathbb X$, we have
\begin{equation} \label{fct-beta2}
\forall n\geq1,\quad S_n+\chi(X_n)-\chi(x)\in b_n + \SS_0\ \ {\mathbb P}_{(x,0)}-a.s.. 
\end{equation}
Hence, a sufficient condition for $(S_n)_n$ to be non-sublattice in $\SS$ is that there exists no $(a_1,\SS_0,\chi(\cdot))$ 
with $a_1\in{\mathbb R}^2$, $\SS_0$ a proper subgroup of $\SS$, $\chi$ a bounded measurable function  
from $\mathbb X$ to $\mathbb R^2$, satisfying for $\pi$-almost every $x\in\mathbb X$, 
$$ S_1+\chi(X_1)-\chi(x)\in a_1 + \SS_0\ \ {\mathbb P}_{(x,0)}-a.s.. 
$$
In some cases (such as additive functionals, or general MRW with 
$\mathcal B$ in {\bf (A1)}-{\bf (A4)} composed of classes of functions modulo $\pi$), the last condition
is equivalent to the non-sublattice condition, see Remark~\ref{dern-rque}. 

\noindent{\scriptsize $\bullet$} {\it Remarks for Markov additive functionals.}  
Let $(X_n)_{n\in\NN}$ be 
a Markov chain  with state space $\X$, transition kernel $Q$, 
invariant distribution $\pi$, and initial distribution $\mu$. 
Here, 
given $\xi=(\xi_1,\xi_2) : \X\rightarrow \RR^2$ 
a $\pi$-centered  function (i.e.~for $i=1,2$, $\xi_i$ is $\pi$-integrable and $\pi(\xi_i)=0$), we consider the classical MRW $(X_n,S_n)_{n\in\NN}$ 
defined by $S_0=0$ and $\forall n\geq1$: 
\begin{equation} \label{add-fct}
S_n:=\sum_{k=1}^n\xi(X_k).
\end{equation}
The sequence $(S_n)_n$ is called an {\it additive functional} (AF) of $(X_n)_n$. In this case, the two following remarks are of interest.

\begin{rqe1} \label{sec:CLT-xi} (Reduction of {\bf (A6)}.)\\
Condition~(\ref{Sn-supported-by-H}) holds if and only if $\xi(\X)\subset\SS$. 
Under this assumption, $(S_n)_n$ is non-sublattice in $\SS$ if and only if  
there exists no $(a_0,\SS_0,A,\chi)$ with 
$a_0\in\RR^2$, $\SS_0$ a proper closed subgroup in $\SS$,  
$A\in{\mathcal X}$ a $\pi$-full
$Q$-absorbing set (i.e.~such that $\pi(A)=1$ and $Q(z,A)=1$ for all $z\in A$), 
$\chi\, :\, \X\rightarrow\RR^2$ a bounded measurable function, such that 
$$\forall x\in A, \quad \xi(y)+\chi(y)-\chi(x)\in a_0 + \SS_0\ \ Q(x,dy)-a.s..$$ 
This statement has been proved in \cite[Section 5.2]{loicsoaz} 
when $\SS=\RR^2$, extension to Cases {\it (H2)-(H3)} is easy. 
\end{rqe1}
\begin{rqe1} \label{rk-non-matrix}
If $(S_n)_n$ satisfies the standard CLT ($A_n=\sqrt n$) and $(S_n)_n$ is non-sublattice in $\SS$, then the covariance matrix $\Gamma$ of the CLT is automatically positive definite, see e.g.~\cite[Section~5.2]{loicsoaz}.   
\end{rqe1}

%
%
%
%
\section{Proof of Theorems I-II} \label{sec:basic-fact}
%
%
Let $(X_n,S_n)_{n\in \NN}$ be a Markov random walk with state space $\X\times\RR^2$, 
and let $\SS$ be a two-dimensional closed subgroup of $\RR^2$ 
satisfying (\ref{Sn-supported-by-H}). We use the notations of Section~\ref{presentation}.

The first assertion of Theorem~I is established in Subsection~\ref{sub-rec-stato}, the second one in Subsection~\ref{sub-sta-nonstat}. Theorem~II is proved in Subsection~\ref{sub-KochenStone}. Auxiliary statements of interest are also presented in these subsections. 
\subsection{Recurrence set in the stationary case (proof of Theorem I-(a))} \label{sub-rec-stato}
To prove Theorem~I, we define the r.v.~$\xi_0=0$ and $\xi_k = S_{k}-S_{k-1}$ for $k\geq1$. 
{}From the additive property~(\ref{Add_pro}), it can be easily seen that 
the distribution of $((X_{n+k},\xi_{n+k}))_{k\ge 1}$
given $\{X_n=x,\ S_n=s\}$ is equal to the distribution of $((X_{k},\xi_{k}))_{k\ge 1}$ 
under ${\mathbb P}_{(x,0)}$. We assume 
(without loss of generality) that $((X_n,\xi_{n}))_{n\geq0}$ 
is the canonical Markov chain with transition kernel $\tilde{P}((x,s);\cdot) := P((x,0);\cdot)$.
Hence, defining the $\sigma$-algebra ${\mathcal F}_n = \sigma(\xi_k,\, 0\leq k \leq n)$ and writing $\theta$ for the usual shift operator on $\Omega = (\X\times{\mathbb R}^2)^{{\mathbb N}}$, we obtain for every bounded measurable function $F:\Omega\rightarrow\mathbb R$ and for every $x\in\X$: ${\mathbb E}_{(x,0)}[F\circ \theta^n\, \vert\, {\mathcal F}_n\,]
={\mathbb E}_{(X_n,0)}[F]$.  
\begin{rqe1}\label{mark-prop}
For $A\in \B(\RR^2)$, $k\in\NN^*$, set $Y_k = \prod_{j=k}^{+\infty} {\mathbf 1}_A(S_j)$, 
and $f_k(x) = {\mathbb E}_{(x,0)}[Y_k]\ $ ($x\in \X)$. Then, for any $B\in \B(\RR^2)$ and 
$n\in{\mathbb N}^*$, we have 
$${\mathbb P}_{(\pi,0)}\bigg(S_n\in B,\ \ S_{n+j}-S_n\in A,\ \forall j\geq k\bigg) = 
{\mathbb E}_{(\pi,0)}\big[{\mathbf 1}_B(S_n)\, f_k(X_n)\big].$$ 
\end{rqe1}
\noindent Note that, for any $A\in \B(\RR^2)$, the corresponding function $f_k$ in 
Remark~\ref{mark-prop} is nonnegative, bounded and measurable. We start by proving 
the recurrence of $(S_n)_n$. 
\begin{lem1} \label{0-in-R}
We have: $0\in{\mathcal R}_{(\pi,0)}$.
\end{lem1}
\begin{proof} 
Let $\varepsilon >0$, $k\geq 1$. 
Let us prove 
that ${\mathbb P}_\pi(\, \exists j\geq k\, :\, |S_j| <2\varepsilon\, ) = 1$. 
For any $n\geq 1$, set 
$$A^{(k)}_n = \bigg[\, |S_n| <\varepsilon,\ |S_{n+j}| \geq \varepsilon,\ 
\forall j\geq k\, \bigg].$$ 
If $|n-n'| \geq k$, then $A^{(k)}_n\cap A^{(k)}_{n'}=\emptyset$. Hence we have 
$\sum_{n\geq1}{\mathbb P}_{(\pi,0)}(A^{(k)}_n) \leq k$. Moreover we have 
$${\mathbb P}_{(\pi,0)}(A^{(k)}_n) \geq {\mathbb P}_{(\pi,0)}\bigg(\, |S_n| <\varepsilon,\ |S_{n+j}-S_n| \geq 2\varepsilon,\ \forall j\geq k\, \bigg).$$ 
Then, applying Remark~\ref{mark-prop}  with $B=\{z\in\RR^2 : |z| < \varepsilon\}$ and $A = \{z\in\RR^2 : |z| \geq 2\varepsilon\}$, 
we obtain $\sum_{n\geq1}{\mathbb E}_{(\pi,0)}\big[{\mathbf 1}_B(S_n)\, f_k(X_n)\big]\leq k$. 
But (\ref{LLT-f-bornee}) gives as $n\to+\infty$: 
$${\mathbb E}_{(\pi,0)}[{\mathbf 1}_B(S_n)\, f_k(X_n)]\sim a_n\, \pi(f_k) m_\SS(B).$$  
Since $0\in \SS$ and $B$ is centered at 0, we have $m_\SS(B) >0$. Finally the fact that $\sum_{n\ge 1}a_n=\infty$ implies $\pi(f_k) = {\mathbb P}_{(\pi,0)}(Y_k=1) =0$.
\end{proof}
\begin{proof}[\bf Proof of assertion {\bf (a)} in Theorem~I]
Let $s\in \SS$. Let us show that 
$$\forall \varepsilon >0,\ \forall k\geq 1,\ \ {\mathbb P}_{(\pi,0)}\big(\, |S_j-s| \geq 2\varepsilon,
\ \forall j\geq k\, \big) = 0.$$
Let $\varepsilon >0$ and $k\geq 1$ be fixed. Set $B'=\{z\in\RR^2 : |z+s| < \varepsilon\}$ and $A' = \{z\in\RR^2 : |z-s| \geq 2\varepsilon\}$, and denote by $Y_k'$ and $f_k'$ 
the elements associated to $A'$ as in Remark~\ref{mark-prop}. Then,
according to Lemma \ref{0-in-R} and Remark~\ref{mark-prop}, we have for $n\geq 1$ 
\begin{eqnarray*}
0={\mathbb P}_{(\pi,0)}\big(\, |S_{n+j}| \geq \varepsilon,\ \forall j\geq k\, \big) &\geq& 
{\mathbb P}_{(\pi,0)}\big(\, S_n\in B',\, S_{n+j}- S_n\in A',\, \forall j\geq k\, \big) \\ 
&=& \ {\mathbb E}_{(\pi,0)}\big[{\mathbf 1}_{B'}(S_n)\, f_k'(X_n)\big]. 
\end{eqnarray*}
Hence ${\mathbb E}_{(\pi,0)}\big[{\mathbf 1}_{B'}(S_n)\, f_k'(X_n)\big]=0$. From (\ref{LLT-f-bornee}) and $m_\SS(B')>0$ (since $B'$ is centered at $-s\in \SS$), it then follows that $\pi(f'_k) = {\mathbb P}_\pi(Y_k'=1) = 0$. 
\end{proof} 
\subsection{From stationarity to non-stationarity under Harris recurrence
(proof of Theorem I-(b))} \label{sub-sta-nonstat} 
Let us define the following subset of $\X$:  
$${\mathcal A} := \{x\in\X : {\mathcal R}_{(x,0)} = \SS\}.$$ 
Property ${\mathcal R}_{(\pi,0)} = \SS$ implies that $\pi({\mathcal A}) = 1$
(since $\SS$ is separable).
Of course, if ${\mathcal A}=\X$, we obtain ${\mathcal R}_{(\mu,0)}=\SS$ 
for every initial distribution $\mu$ of the driving Markov chain $(X_n)_{n\in\NN}$. 

The second assertion of Theorem~I follows from the following statement. 
\begin{prop1} \label{pro-sta-nonstat} The following assertions hold: 
\begin{enumerate}
\item[(i)] If $x\in\X$ is such that ${\mathbb P}_{(x,0)}\big(X_n\in {\mathcal A}\ \text{i.o.}\big)=1$, 
then $x\in{\mathcal A}$ (i.e. ${\mathcal R}_{(x,0)}=\SS$). 
\item[(ii)] If the driving Markov chain $(X_n)_{n\in\NN}$ is Harris recurrent and  
if ${\mathcal R}_{(\pi,0)}=\SS$, then ${\mathcal A}=\X$. In this case, we have  
${\mathcal R}_{(\mu,0)}=\SS$ for every initial distribution $\mu$ of $(X_n)_{n\in\NN}$. 
\end{enumerate}
\end{prop1}
\begin{lem1} \label{doob-additive}
Let $(s,\varepsilon)\in \SS\times(0,1)$ and $E := [|S_n-s| < \varepsilon\ \text{i.o.}]$. 
For every $x\in \X$, we have for ${\mathbb P}_{(x,0)}-$almost every $\omega\in\Omega$: 
$$\lim_{k\rightarrow +\infty} {\mathbb P}_{(X_k(\omega),0)}\bigg(\big\vert S_n-
\big(s-S_k(\omega)\big)\big\vert
<\varepsilon\ \text{i.o.}\bigg) = {\bf 1}_{E}(\omega). $$
\end{lem1}
\begin{proof}[Proof of Lemma \ref{doob-additive}]
Let $x\in\mathbb X$.
According to a classical argument due to Doob (see \cite[Prop. V-2.4]{Neveu}), we have for ${\mathbb P}_{(x,0)}-$almost every $\omega\in\Omega$: 
$\lim_{k\rightarrow +\infty} {\mathbb P}_{(X_k(\omega),S_k(\omega))}(E)={\bf 1}_{E}(\omega)$. Then the desired property easily follows from 
the additive property (\ref{Add_pro}). 
\end{proof}
\begin{proof}
[Proof of Proposition~\ref{pro-sta-nonstat}]
We suppose that $(X_n,S_n)_n$ is the canonical version defined on the set 
$\Omega:=(\X\times\RR^2)^{\NN}$. Let us fix any $(s,\varepsilon)\in \SS\times(0,1)$, and set 
$E := \big[|S_n-s| < \varepsilon\ \text{i.o.}\big]$. 
Using  the assumption in {\it (i)}, Lemma~\ref{doob-additive} and Lebesgue's theorem,
using finally the definition of ${\mathcal A}$ and the fact that  $S_k-s\in \SS$ $\, {\mathbb P}_{(x,0)}-$a.s. (use (\ref{Sn-supported-by-H})), we obtain the following property  
$${\mathbb P}_{(x,0)}(E) = \lim_k \int_{\{\omega : X_n(\omega)\in {\mathcal A}\ \text{i.o.}\}} {\mathbb P}_{(X_k(\omega),0)}\left( \big\vert S_n-\big(s-S_k(\omega)\big)\big\vert
<\varepsilon\ \text{i.o.}\right)\, d{\mathbb P}_{(x,0)}(\omega) = 1,$$
from which we deduce ${\mathcal R}_{(x,0)} = \SS$. Now, if ${\mathcal R}_{(\pi,0)} = \SS$, then $\pi({\mathcal A}) = 1$, so that the Harris recurrence of $(X_n)_{n\in\NN}$ gives ${\mathbb P}_{(x,0)}\big(X_k\in {\mathcal A}\ \text{i.o.}\big) = 1$ for all $x\in\X$. Thus (ii) follows from (i). 
\end{proof}
Lemma~\ref{doob-additive}, based on both Markov and additive properties of $(X_n,S_n)_n$, plays an important role in the previous proof, as well as in the main statement (Proposition~\ref{rec-1}) of the next section. 
\subsection{Borel-Cantelli adaptation of Kochen and Stone (Proof of Theorem~II)} \label{sub-KochenStone}
We present now a general strategy to obtain ${\mathcal R}_{(x,0)}=\SS$ for every 
$x\in\mathbb X$, even 
when the driving Markov chain is not Harris-recurrent. In particular Theorem~II 
directly follows from the next Corollary~\ref{coro-rqe-cond-Koc-sto} and 
Proposition~\ref{rec-1}. 

The following Proposition~\ref{prop-KochenStone}, as well as its 
Corollary~\ref{coro-rqe-cond-Koc-sto}, are true for any sequence $(Y_n)_{n\in\NN}$ of r.v.~defined on some probability space 
$(\Omega, {\mathcal F},{\mathbb P})$ and taking their values in $\RR^2$.

\begin{prop1}[\cite{KochenStone}] \label{prop-KochenStone}
Let $(s,\varepsilon)\in{\mathbb R}^2\times(0,1]$. Assume that there exists 
$c\in[1,+\infty)$ such that
\begin{subequations}
\begin{eqnarray}
& & \ \ \sum_{n\ge 1}{\mathbb P}\big(\vert Y_n-s\vert<\varepsilon\big) = \infty 
\label{cond-KochenStone-1}\\
& & \ \ \liminf_{N\rightarrow +\infty}
\frac{ \sum_{n,m=1}^N{\mathbb P}\big( \vert Y_n-s\vert<\varepsilon,\  \vert Y_{m}-s\vert
<\varepsilon\big)}
{\left(\sum_{n=1}^N{\mathbb P}\big( \vert Y_n-s\vert<\varepsilon\big)  \right)^2}
\le c. \label{cond-KochenStone-2}
\end{eqnarray}
\end{subequations}
Then we have: $\ \displaystyle {\mathbb P}\big(|Y_n-s| < \varepsilon\ \text{i.o.}\big)\ge
 \frac{1}{c}$. 
\end{prop1}
\begin{coro1} \label{coro-rqe-cond-Koc-sto} 
Let $(s,\varepsilon)\in{\mathbb R}^2\times(0,1]$. Assume that Condition~(\ref{cond-KochenStone-1}) is fulfilled and that there exists $d\in(0,+\infty)$ such that
\begin{equation} \label{rqe-cond-Koc-sto}
\liminf_{N\rightarrow +\infty}
\frac{ \sum_{n,m=1}^N{\mathbb P}\big( \vert Y_n-s\vert<\varepsilon,\, \vert Y_{n+m}-s\vert
<\varepsilon\big)}
{\left(\sum_{n=1}^N {\mathbb P}( \vert Y_n-s\vert<\varepsilon)  \right)^2}
\le d. 
\end{equation}
Then we have: $\ \displaystyle {\mathbb P}\big(|Y_n-s| < \varepsilon\ \text{i.o.}\big)\ge \frac{1}{2d}$. 
\end{coro1}
\begin{proof}[Proof of Corollary~\ref{coro-rqe-cond-Koc-sto}] 
Let us define $p_{n,m}:={\mathbb P}( \vert Y_n-s\vert<\varepsilon,\  \vert Y_{n+m}-s\vert
<\varepsilon)$. Observe that $p_{n,0} := {\mathbb P}( \vert Y_n-s\vert<\varepsilon)$. We have 
$$\sum_{n,m=1}^N {\mathbb P}(|Y_n-s|<\varepsilon,\
  \vert Y_m -s\vert<\varepsilon) \leq 2 \sum_{n=1}^N\sum_{m=n}^N 
   p_{n,m-n}\le 2\sum_{n=1}^N\sum_{m=0}^N 
   p_{n,m} = 2\left(\sum_{n,m=1}^N p_{n,m} +  \sum_{n=1}^N    p_{n,0}\right).$$
{}From (\ref{cond-KochenStone-1}) and the previous inequality, we obtain (\ref{cond-KochenStone-2}) with $c=2d$. 
\end{proof}
Notice that, even in the i.i.d.~case,
Corollary~\ref{coro-rqe-cond-Koc-sto} does not give 
$p_n:={\mathbb P}(|Y_n-s| < \varepsilon\ \text{i.o.}) = 1$ as expected, but only 
$p_n\ge 1/2$ (since constant $d$ is equal to 1). 
Therefore, further arguments (here based on the additive property (\ref{Add_pro})) must be exploited to deduce 
the recurrence set from Corollary~\ref{coro-rqe-cond-Koc-sto}. 
The next proposition gives such 
a result for general Markov random walks. 

Again $(X_n,S_n)_{n\in \NN}$ denotes a MRW with state space $\X\times\RR^2$, and $\SS$ is given in (\ref{Sn-supported-by-H}). 
\begin{prop1} \label{rec-1}
Let $\varepsilon>0$. Assume that there exists a real number $e_\varepsilon>0$ such that
\begin{equation} \label{min-c-sur-A}
\forall (x',s')\in \X\times \SS,\ \ 
{\mathbb P}_{(x',0)}\big(|S_n-s'| < \varepsilon\ \text{i.o.}\big)\ge e_\varepsilon.
\end{equation}
Then, for every $(x,s)\in \X\times \SS$, we have: ${\mathbb P}_{(x,0)}\big(|S_n-s| < \varepsilon\ \text{i.o.}\big) = 1$. 
In particular, if (\ref{min-c-sur-A}) is fulfilled for every $\varepsilon\in(0;1)$, 
then we have for every $x\in\X$: $\, {\mathcal R}_{(x,0)} = \SS$ . 
\end{prop1}
\begin{proof}
Suppose that $(X_n,S_n)_n$ is the canonical version defined on 
$\Omega:=(\X\times\RR^2)^{\NN}$. Let us fix $(x,s)\in \X\times \SS$, and set $S := [|S_n-s| < \varepsilon\ \text{i.o.}]$. 
Then, from Lemma~\ref{doob-additive}, (\ref{Sn-supported-by-H}) and (\ref{min-c-sur-A}), 
it follows that, for ${\mathbb P}_{(x,0)}-$almost every $\omega\in\Omega$, 
we have: ${\bf 1}_{S}(\omega)\geq e_\varepsilon$. Hence: ${\bf 1}_{S}=1$ 
$\, {\mathbb P}_{(x,0)}$-a.s.. 
\end{proof}
%
\section{LLTs under Hypotheses~{{\bf (A1)}}-{{\bf (A6)}} (Proof of Theorems~III-IV)} \label{sec-meth-spec}
%
Let $(X_n,S_n)_{n\in \NN}$ be a Markov random walk with state space $\X\times\RR^2$, 
let $\SS$ be a two-dimensional closed subgroup of $\RR^2$ 
satisfying (\ref{Sn-supported-by-H}). 
Hypotheses of Theorems~I-II involve some local limit properties. This is obvious for Theorem~I since 
Hypothesis~(\ref{LLT-f-bornee}) directly writes as a local limit property. The next lemma shows that this is also true for Theorem~II, more precisely: Conditions~(\ref{cond-cor-resume-1}) and (\ref{cond-cor-resume-2}) of Section~\ref{presentation} are implied by the limit properties (\ref{LLT-1-kochen}) and (\ref{LLT-2-kochen}). 

The notations $u_n\sim v_n$ or $u_n\sim_n v_n$ refer to the usual equivalence relation between two sequences as $n\rightarrow +\infty$. 
We write $u_{n,m} \sim_{(n,m)} v_{n,m}$ when, for $n$, $m$ large enough, we have $u_{n,m} = v_{n,m}(1+\eta_{n,m})$ for some bounded $(\eta_{n,m})_{n,m}$ such that $\eta_{n,m} \rightarrow 0$ when 
$\min(n,m)\rightarrow+\infty$.

\begin{lem1} \label{sim-sum}
Let $\mu$ be a probability measure on $\X$ and let $B$ be a ball in $\RR^2$. 
Assume that there exist a constant $D>0$, a sequence $(a_n)_{n\ge 1}$ of positive numbers such that $\sum_{n\ge 1}a_n=\infty$ and:  
\begin{subequations}
\begin{eqnarray}
& & \ \ {\mathbb P}_{(\mu,0)}\big(S_n\in B\big) \sim_n D\, a_n\,  m_\SS(B) \label{equi-1} \\ 
& & \ \ {\mathbb P}_{(\mu,0)}\big((S_n,S_{n+m})\in B^2\big) \sim_{(m,n)}  D^2\, a_na_m\,  m_\SS(B)^2. \label{equi-2} 
\end{eqnarray}
\end{subequations}
Then we have 
$$\sum_{n\ge 1}{\mathbb P}_{(\mu,0)}\big(S_n\in B\big) = \infty, $$
$$\lim_{N\rightarrow +\infty}
 \frac{\sum_{n,m=1}^{N} {\mathbb P}_{(\mu,0)}\big((S_n,S_{n+m})\in B^2\big)}{\left(\sum_{n=1}^{N} {\mathbb P}_{(\mu,0)}\big(S_n\in B\big)\right)^{2}} = 1. $$
\end{lem1}
\begin{proof} 
Set $p_{n,m} := {\mathbb P}_{(\mu,0)}((S_n,S_{n+m})\in B^2)$. Note that $p_{n,0} := {\mathbb P}_{(\mu,0)}\big(S_n\in B\big)$.  We have when $N\rightarrow+\infty$: 
$$\sum_{n=1}^N p_{n,0} \sim_N D\, m_\SS(B)\,  \sum_{n=1}^N  a_n\ \ \ \ \ \mbox{and}\ \ \ \ 
\sum_{n,m=1}^N p_{n,m} \sim_N D^2\, m_\SS(B)^2\, \left(\sum_{n=1}^N a_n\right)
      \left(\sum_{m=1}^N a_m\right), $$
from which we deduce the desired statement.  
\end{proof}

In the next Propositions~\ref{sec:loi_stable-ter}-\ref{sec:loi_stable}, the local limit properties (\ref{LLT-f-bornee}) and (\ref{equi-1})-(\ref{equi-2}) are obtained under 
Hypotheses~{{\bf (A1)}}-{{\bf (A6)}}. Theorems~III-IV are then deduced from Theorems~I-II. Another interesting application to recurrence is presented in Corollary~\ref{rec-spectral-1}.  

Recall that $B:= B(s,\varepsilon)$ is the open ball in $\RR^2$, centered at $s$ with radius $\varepsilon$. 
The sequence $(A_n)_n$ is given in {\bf (A5)}. Let us define the following positive constant: 
$D_\SS := (2\pi)^{-1}\, c_\SS\, (\det\Gamma)^{-1/2}$, where $c_\SS=\varepsilon_\SS$ for (H1) (H2), and $c_\SS = |\det B|$ for (H3), where Cases (H1) (H2) (H3) are described in Section~\ref{simp-rks-S}. 
\begin{prop1}\label{sec:loi_stable-ter} 
Assume that Hypotheses~{{\bf (A1)}}-{{\bf (A6)}} hold true and that $\mu\in(\widehat{\mathcal B})'$. Then, for every $(s,\varepsilon)\in \SS\times(0;\varepsilon_\SS)$, for every bounded
nonnegative $f\in\mathcal B$, we have: 
\begin{equation} \label{LLT-f-bornee-gene}
\EE_{(\mu,0)}\big[f(X_n)\, {\mathbf 1}_{B(s,\varepsilon)}(S_n)\big] \sim_n D_\SS\, \pi(f)\, A_n^{-2}\, m_\SS(B(s,\varepsilon)).
\end{equation}
\end{prop1}
\begin{prop1} \label{sec:loi_stable} 
Assume that  Hypotheses~{{\bf (A1)}}-{{\bf (A6)}} hold true and that $\mu\in(\widehat{\mathcal B})'$. 
Then, for every $(s,\varepsilon)\in \SS\times(0;\varepsilon_\SS)$, we have (\ref{equi-1})-(\ref{equi-2}) with $B:=B(s,\varepsilon)$, $a_n:=A_n^{-2}$ and $D:=D_\SS$. 
\end{prop1}

Propositions~\ref{sec:loi_stable-ter} and \ref{sec:loi_stable} are proved in Section~\ref{sec-proof-deux-prop}. 

\begin{proof}[Proof of Theorem~III-IV]
Condition~(\ref{LLT-f-bornee}) of Theorem~I is nothing else but (\ref{LLT-f-bornee-gene}) stated with $\mu=\pi$ and for every nonnegative bounded measurable function $f : \X\rightarrow \RR$. Note that $\pi\in(\widehat{\mathcal B})'$ from the continuous inclusion $\widehat{\mathcal B} \subset {\mathbb L}^1(\pi)$. Consequently 
Theorem~III follows from Proposition~\ref{sec:loi_stable-ter}. Theorem~IV follows from 
Proposition~\ref{sec:loi_stable}, Lemma~\ref{sim-sum} and Theorem~II.
\end{proof}

When the assumption of Theorem~IV on Dirac distributions is not fulfilled, the following corollary may also be of interest. It follows from 
Proposition~\ref{sec:loi_stable}, Lemma~\ref{sim-sum} and Corollary~\ref{coro-rqe-cond-Koc-sto}. 
\begin{coro1} \label{rec-spectral-1}
Assume that Hypotheses~{{\bf (A1)}}-{{\bf (A6)}} hold true, that $\mu\in(\widehat{\mathcal B})'$. Then we have for every $(s,\varepsilon)\in \SS\times (0;\varepsilon_\SS]$: 
$\displaystyle \ {\mathbb P}_{(\mu,0)}\big(|S_n-s| < \varepsilon\ \text{i.o.}\big) \geq 1/2$. 
\end{coro1}

We present now some remarks concerning the operator-type Hypotheses {\bf (A1)}-{\bf (A4)} of Section~\ref{presentation}. Further comments can be found in \cite[Sect.~4-5]{loicsoaz}. Actually Hypotheses~{\bf (A1)} and {\bf (A2)}-{\bf (A3)} are the key assumptions of the weak 
Nagaev-Guivarc'h spectral method presented in \cite[Sect.~4-5]{loicsoaz}, which is used in Section~\ref{sec-proof-deux-prop} to prove  
Propositions~\ref{sec:loi_stable-ter} and \ref{sec:loi_stable}. 

\noindent {\scriptsize $\bullet$} {\bf Comments on Hypotheses~{\bf (A1)}-{\bf (A4)}.} The strong ergodicity condition~{\bf (A1)} only involves the driving Markov chain $(X_n)_{n\in\NN}$ of the MRW. 
More specifically, defining the following rank-one projection in ${\mathcal L}({\mathcal B})$,  
\begin{equation} \label{Pi}
\forall f\in{\mathcal B},\quad \Pi f = \pi(f){\mathbf 1}_{\mathbb X},  
\end{equation}
Condition~{\bf (A1)} writes as: $\lim_n\|Q^n-\Pi\|_{{\mathcal B}} = 0$. 
This can be easily seen that the last condition is equivalent to 
$\|Q^n-\Pi\|_{{\mathcal B}} = O(\kappa^{n})$ for some $\kappa\in(0,1)$. 
Under Condition~{\bf (A1)}, the technical condition {\bf (A4)} is satisfied in many cases, see \cite[p.~436]{loicsoaz}. Mention that Inequalities $|\lambda|^n|f|\le Q^n|f|$ and $\vert f\vert\le\pi(|f|)$ in {\bf (A4)} must be understood as follows: they hold, either everywhere on $\X$ if ${\mathcal B}$ is a space of functions, or $\pi$-almost everywhere on $\X$ if ${\mathcal B}$ is a space of classes modulo $\pi$. 

\noindent The condition $Q(t)\in{\mathcal L}({\mathcal Y})$ for ${\mathcal Y}={\mathcal B}$ or ${\mathcal Y}=\widehat{\mathcal B}$ means that, for every $f\in \mathcal Y$, 
function  $\X\ni x\mapsto \big(Q(t)f\big)(x)$ 
(or its class mod.~$\pi$) belongs to $\mathcal Y$,
and that $f\mapsto Q(t)f$ is in ${\mathcal L}({\mathcal Y})$. 

\noindent Hypotheses~{\bf (A2)}-{\bf (A3)} enable the use of the Keller-Liverani perturbation theorem in the Nagaev-Guivarc'h spectral method. 
Note that Hypotheses~{\bf (A2)}-{\bf (A3)} involve not only the transition kernel $Q$ of the driving Markov chain $(X_n)_{n\in\NN}$, but also the additive component $S_n$ of the MRW. For instance, if $S_n:=\sum_{k=1}^n\xi(X_k)$ is an additive functional, then 
{\bf (A2)}-{\bf (A3)} mainly focus on the function $\xi : \X\rightarrow\RR^2$. 

\noindent Hypothesis~{\bf (A2)} is a continuity condition involving two different spaces ${\mathcal B}\subset\widehat{\mathcal B}$. This condition is much weaker than in the usual perturbation theorem involving a single space (i.e.~$\widehat{\mathcal B}={\mathcal B}$): for instance, as illustrated in \cite[Sect.~3]{loicsoaz}, Hypothesis~{\bf (A2)} does not hold in general with $\widehat{\mathcal B}={\mathcal B}$
in the classical Markov models considered in Applications~\ref{aplli-rho-mix}-\ref{appli-non-stand} of Section~\ref{presentation}. Condition~(\ref{ineg-DF}) in {\bf (A3)} is the so-called Doeblin-Fortet inequality: here it is required for all the $Q(t)$ in a uniform way on compact sets of $\RR^2$.  

\noindent Finally, concerning the notion of essential spectral radius, recall that $T\in{\mathcal L}({\mathcal B})$ is said to be quasi-compact if there exist $r_0\in(0,1)$, $m\in\NN^*$, $\lambda_i\in\CC$, $p_i\in\NN^*$ ($\, i=1,\ldots,m$) such that: 
$$\mathcal B = \overset{m}{\underset{i=1}{\oplus}} \ker(T - \lambda_i I)^{p_i}\, \oplus H,$$
where 
$|\lambda_i| \geq r_0$, $\, 1\le \dim\ker(T-\lambda_i I)^{p_i} < \infty$, and $H$ is a closed $T$-invariant subspace such that $\sup_{h\in H,\, \|h\|_{\mathcal B}\leq1}\|T^nh\|_{\mathcal B} = O({r_0}^n)$.
If $T$ is quasi-compact, then $r_{ess}(T)$ is the infimum bound of the real numbers $r_0$ such that the last conditions hold. If $T$ is not quasi-compact, then $r_{ess}(T)$ is equal to the spectral radius of $T$. For further details on the essential spectral radius, in particular for the link with the Doeblin-Fortet inequalities, see \cite{hen2,hulo}.

\section{Proof of Propositions~\ref{sec:loi_stable-ter} and \ref{sec:loi_stable}} \label{sec-proof-deux-prop} 
\subsection{Spectral properties of $Q(t)$ under Hypotheses~{{\bf (A1)}}-{{\bf (A6)}}} \label{subsub-spectal} 
The Fourier kernels $Q(t)$ are defined in (\ref{Fourier}). Note that $Q(0) = Q$, and $Q(t+g)=Q(t)$ for every $(t,g)\in {\mathbb R}^2\times \SS^*$. The positive definite symmetric $2\times2$-matrix $\Gamma$ and the slowly varying function $L(\cdot)$ used below are defined in {\bf (A5)}. Recall that we have set in (\ref{Pi}): $\forall f\in{\mathcal B},\ \Pi f = \pi(f){\mathbf 1}_{\mathbb X}$. 
\begin{prop1} \label{sec:KL} 
Under Hypotheses~{{\bf (A1)}}-{{\bf (A6)}}, the following assertions hold true:
\begin{enumerate}[(a)]
\item There exist two real numbers $\alpha>0$ and 
$\kappa\in[0,1)$, a function $t\mapsto\lambda(t)$ from $B(0,\alpha)$ into $\mathbb C$, 
a bounded map $t\mapsto \Pi(t)$ from $B(0,\alpha)$ into ${\mathcal L}({\mathcal B})$ 
such that $\lambda(0)=1$, $\Pi(0) = \Pi$ and 
\begin{equation} \label{dec-KL}
\sup_{ t\in B(0,\alpha)}
   \left\Vert Q(t)^n-\lambda(t)^n\Pi(t)\right\Vert_{\mathcal B}=O(\kappa^n),
\end{equation}
\begin{equation} \label{Dl-2}
\lambda(t)=1- \frac{1}{2}\langle t,\Gamma t\rangle L(\vert t\vert^{-1})\big(1+\varepsilon(t)\big),
\end{equation}
\item If $\mu\in(\widehat{\mathcal B})'$, then  
\begin{equation} \label{b-stat-1}
\forall f\in\mathcal B,\ \ \lim_{t\rightarrow0}\, \mu\big(\Pi(t)f\big) = \pi(f) 
\end{equation}
\begin{equation} \label{b-stat} 
\lim_{(u,v)\rightarrow0}\, \mu\big(\Pi(u)\Pi(v){\bf 1}_\X\big) =  1. 
\end{equation}
\item For any compact subset $K$ of $\RR^2\setminus \SS^*$, there exists $\rho=\rho(K)\in[0,1)$ such that 
\begin{equation} \label{non-ari-spectral-lattice} 
\sup_{t\in K}\Vert Q(t)^n\Vert_{\mathcal B} = O(\rho^n).
\end{equation}
\end{enumerate}
\end{prop1}

\begin{proof}
\noindent Property~(\ref{dec-KL}) is presented in \cite[Sect.~4]{loicsoaz}. 
Property~(\ref{Dl-2}) is established in \cite[lem.~4.2]{ihp1} when the standard CLT holds in {\bf (A5)}, see also \cite[Lem.~5.2]{loicsoaz}. 
Extension to non-standard CLT is easy, see Appendix~\ref{conv-gauss-mambda}. 
Property~(\ref{non-ari-spectral-lattice}) is established in \cite[Sect.~5]{loicsoaz} when $\SS=\RR^2$. Extension to a proper closed subgroup $\SS$ of $\RR^2$ is simple, see Appendix~\ref{app-B}. To obtain (\ref{b-stat-1})-(\ref{b-stat}), recall that the main argument in the proof of (\ref{dec-KL}) is the Keller-Liverani perturbation theorem  \cite{KL}, which also gives the following properties (see  \cite[Sect.~4]{loicsoaz} for details): 
$$M := \sup_{ t\in B(0,\alpha)}\|\Pi(t)\|_{\mathcal B} < \infty\quad \text{ and } \quad \lim_{t\rightarrow 0}\|\Pi(t)-\Pi\|_{{\mathcal B},\widehat{\mathcal B}} =0.$$  
Then (\ref{b-stat-1}) follows from $\mu\in(\widehat{\mathcal B})'$ and the last property. In particular, since  $\pi\in(\widehat{\mathcal B})'$, (\ref{b-stat-1}) holds with $\mu=\pi$. Next, using $\Pi(\Pi(v) {\bf 1}_\X) = \pi(\Pi(v) {\bf 1}_\X){\bf 1}_\X$, $\Pi{\bf 1}_\X = {\bf 1}_\X$ and $\mu({\bf 1}_\X)=1$, we obtain:  
\begin{eqnarray*}
 \big|\mu\big(\Pi(u)\Pi(v){\bf 1}_\X\big) - 1 \big| 
 &\leq& \mu\big(\, \big| \big(\Pi(u)-\Pi\big)\Pi(v){\bf 1}_\X\big|\, \big) +   \big|\pi\big(\Pi(v){\bf 1}_\X\big)-1\big| \\
&\leq& \|\mu\|_{\widehat{\mathcal B}}\ \big\|\Pi(u)-\Pi\big\|_{{\mathcal B},\widehat{\mathcal B}}\ 
\big\|\Pi(v){\bf 1}_\X\big\|_{{\mathcal B}} + \big|\pi\big(\Pi(v){\bf 1}_\X\big)-1\big| \\
&\leq& \|\mu\|_{\widehat{\mathcal B}}\ \big\|\Pi(u)-\Pi\big\|_{{\mathcal B},\widehat{\mathcal B}}\ 
M\, \|{\bf 1}_\X\|_{{\mathcal B}} + \big|\pi\big(\Pi(v){\bf 1}_\X\big)-1\big|. 
\end{eqnarray*}
Hence we have (\ref{b-stat}). 
\end{proof}
\subsection{Preliminary lemmas} \label{subsec-preli-lems} 
Let $f$ be a $\CC$-valued bounded measurable function on $\X$. 
\begin{lem1} \label{fourier}
We have for every $x\in\X$, every $(u,v)\in\RR^2$ and every $(n,m)\in\NN^2$: 
$$\EE_{(x,0)}\big[e^{i\langle u,S_n\rangle}\, e^{i\langle v,S_{n+m}-S_n\rangle}\, 
f(X_{n+m})\big] = \big(Q(u)^n Q(v)^m f\big)(x).$$
Consequently, we have for any initial distribution $\mu$ on $\X$:
$$\EE_{(\mu,0)}\big[e^{i\langle u,S_n\rangle}\, e^{i\langle v,S_{n+m}-S_n\rangle}\, 
f(X_{n+m})\big] = \mu\big(Q(u)^n Q(v)^m f\big).$$
\end{lem1}
\begin{proof}[Proof of Lemma \ref{fourier}]
Using additivity property (\ref{Add_pro}) (see Subsection~\ref{sub-rec-stato}), we obtain 
\begin{eqnarray}
\EE_{(x,0)}\big[e^{i\langle u,S_n\rangle}\, e^{i\langle v,S_{n+m}-S_n\rangle}\, f(X_{n+m})\big] 
&=& \EE_{(x,0)}\bigg[e^{i\langle u,S_n\rangle}\,  
\EE_{(x,0)}\big[e^{i\langle v,S_{n+m}-S_n\rangle}\, f(X_{n+m})\, \big|\, {\mathcal F}_n \big]\bigg] \nonumber  \\
&=& \EE_{(x,0)}\bigg[e^{i\langle u,S_n\rangle}\, \EE_{(X_n,0)}\big[f(X_m)\, 
e^{i \langle v , S_m \rangle }\big]\bigg].  \label{form_map_fourier}
\end{eqnarray}
Applying (\ref{form_map_fourier}) with  $m=1$ and $u=v$, and according to
definition (\ref{Fourier}) of Fourier maps, we obtain for every $n\geq0$,
$$\EE_{(x,0)}\big[e^{i\langle v,S_{n+1}\rangle}\, f(X_{n+1})\big] = \EE_{(x,0)}\bigg[e^{i\langle v,S_n\rangle}\, \EE_{(X_n,0)}\big[f(X_1)\, 
e^{i \langle v , S_1 \rangle }\big]\bigg] = \EE_{(x,0)}\big[e^{i\langle v,S_n\rangle}\big(Q(v)f\big)(X_n)\big].$$
We deduce by induction that we have for all $v\in\RR^2$, $k\geq1$, and for all $\CC$-valued bounded measurable function $g$ on $\X$: 
$$\EE_{(x,0)}\big[e^{i\langle v,S_{k}\rangle}\, g(X_{k})\big] = \big(Q(v)^k g\big)(x).$$
Next, by applying (\ref{form_map_fourier}) (with  $m\geq1$ and $u,v\in\RR^2$) and using the previous equality (first with $g=f$, second with $g=Q(v)^mf$), we obtain 
$$\EE_{(x,0)}\big[e^{i\langle u,S_n\rangle}\, e^{i\langle v,S_{n+m}-S_n\rangle}\, f(X_{n+m})\big] =  
\EE_{(x,0)}\big[e^{i\langle u,S_n\rangle}\, \big(Q(v)^mf\big)(X_n)\big] = \big(Q(u)^nQ(v)^mf\big)(x).$$
\end{proof}
For any Lebesgue-integrable function $h:{\mathbb R}^2\rightarrow\mathbb C$, we define
its Fourier transform $\hat h$ by $\hat h(u):=\int_{{\mathbb R}^2}h(t)e^{-i\langle t,u \rangle}\, dt$, and we set 
$$\forall t\in\RR^2,\ \ P_h(t) := \sum_{g\in \SS^*}\hat h(t+g).$$ 
Let $\mathcal D$ be the fundamental domain of ${\mathbb R}^2/\SS^*$, namely: 
\begin{enumerate}
\item[-] $\ {\mathcal D}:=\RR^2$ in Case~{\it (H1)}, 
\item[-] $\ {\mathcal D}:= [-\frac{a}{2},\frac{a}{2}]\times\RR$ in Case~{\it (H2)}, with $a={2\pi}/b$, 
\item[-] $\ {\mathcal D}:=A([-\frac{1}{2},\frac{1}{2}]^2)$ in Case~{\it (H3)}, with $A:=2\pi(B^*)^{-1}$. 
\end{enumerate}
\begin{lem1}\label{fourier-inverse}
Let $h_1$ and $h_2$ be $\CC$-valued Lebesgue-integrable functions on $\RR^2$ such that their 
Fourier transforms are Lebesgue-integrable on $\RR^2$. 
Then we have for any probability measure $\mu$  on $\X$ and for every $(n,m)\in\NN^2$: 
$$\EE_{(\mu,0)}\big[h_1(S_n)\, f(X_n)\big] = \frac{1}{(2\pi)^{2}}\int_{{\mathcal D}} P_{h_1}(u)\, 
\mu\big(Q(u)^n  f\big)\, du.$$
$$\EE_{(\mu,0)}\big[h_1(S_n)\, h_2(S_{n+m}-S_n)\, f(X_{n+m})\big] = \frac{1}{(2\pi)^{4}}
\int_{{\mathcal D}\times{\mathcal D}}
          P_{h_1}(u)\, P_{h_2}(v)\, \mu\big(Q(u)^n Q(v)^m f\big)\, dudv.$$
\end{lem1}
\begin{proof}[Proof of Lemma \ref{fourier-inverse}]
We easily obtain the first formula by applying the inverse Fourier formula to $h_1$, Lemma~\ref{fourier} (with $m=0$), and finally the fact that $Q(\cdot)$ and $P_{h_1}$ are $\SS^*$-periodic. The second formula can be proved similarly. 
\end{proof}
\begin{lem1} \label{lem-lambda}
Up to a reduction of the positive real number $\alpha$ of Proposition~\ref{sec:KL}, 
there exists $\tilde a>0$ such that, for every $t\in B(0,\alpha)$, we have 
$\vert \lambda(t)\vert\le e^{-\tilde a |t|^2 L(\vert t\vert^{-1})}$, and for all $n\geq1$ 
\begin{equation} \label{lambda-ine-CD}
\left|\lambda\left(\frac{u}{A_n}\right)\right|^n\, {\mathbf 1}_{B(0,\alpha A_n)}(u) 
   \le {\mathbf 1}_{B(0,1)}(u) + e^{-\frac{\tilde a}{4} |u|}\,
    {\mathbf 1}_{\{u\,:\, 1\leq |u| \leq \alpha A_n\}}(u).
\end{equation}
\end{lem1}
\begin{proof}[Proof of Lemma~\ref{lem-lambda}]
{}From (\ref{Dl-2}) and the fact that $\Gamma$ is positive definite, 
there exists $\tilde a>0$ such that, for every $t\in B(0,\alpha)$ (with $\alpha$ possibly reduced), we have 
$$\vert \lambda(t)\vert\le 1 - \tilde a |t|^2 L(\vert t\vert^{-1}) \leq e^{-\tilde 
a |t|^2 L(\vert t\vert^{-1})}.$$
Since $L(\cdot)$ is slowly varying, we know 
(see \cite{Karamata} or \cite{Feller}, p. 282) that there
exist two functions $\ell(\cdot)$ and $\tilde\varepsilon(\cdot)$ such that 
$\lim_{x\rightarrow +\infty}\ell(x)$
exists in $(0,+\infty)$ and $\lim_{x\rightarrow +\infty}\tilde\varepsilon(x)=0$, 
and such that 
\begin{equation}\label{karamata}
L(x)=\ell(x)\exp\left(\int_1^x\frac{\tilde\varepsilon(y)}y\, dy\right). 
\end{equation}
Using this representation of $L$, it is easy to see that
there exists $n_0$ such that, for any $n\ge n_0$ and any $u$ such that
$1\le \vert u\vert\le \alpha A_n$, we have~:
\[
\frac{1}{2}|u|^{-1} \le \frac{L(A_n |u|^{-1})}{L(A_n)}.
\]
{}From $A_n^2\sim nL(A_n)$, one can also assume that, 
for every  $n\ge n_0$ (up to a change of $n_0$), we have $n/A_n^2 \geq \frac{1}{2L(A_n)}$. Therefore we have: $\forall u\in B(0,\alpha A_n)$, $\forall n\geq n_0$, 
\begin{eqnarray*}
\left|\lambda\left(\frac{u}{A_n}\right)\right|^n\, {\mathbf 1}_{B(0,\alpha A_n)}(u) 
&\leq& e^{-n\, \tilde a \frac{|u|^2}{A_n^2} L(A_n|u|^{-1})}\, 
{\mathbf 1}_{B(0,\alpha A_n)}(u)\\ 
&\le& {\mathbf 1}_{B(0,1)}(u) + e^{-\frac{\tilde a}{4} |u|}\, 
{\mathbf 1}_{B(0,\alpha A_n)\setminus B(0,1)}(u).
\end{eqnarray*}
\end{proof}
\subsection{Proof of Proposition~\ref{sec:loi_stable-ter}} \label{proof-sec:loi_stable-ter} 
Let ${\mathcal H}_2$ denote the space of Lebesgue-integrable continuous functions on $\RR^2$ having a compactly supported Fourier transform. Let $f\in\mathcal B$, $f\geq0$ be fixed. Property~(\ref{LLT-f-bornee-gene}) will be proved 
if we establish that we have for every $h\in{\mathcal H}_2$: 
\begin{equation} \label{h-prrof-pn-lattice} 
\lim_n\, 2\pi\, A_n^{2}\, {\mathbb E}_{(\mu,0)}\big[f(X_n)\, h(S_n)\big] =  c_\SS\, (\det\Gamma)^{-1/2}\, \pi(f)\, m_\SS(h). 
\end{equation}
Indeed, (\ref{h-prrof-pn-lattice}) 
ensures that the sequence $(\nu_n)_n$ of positive measures defined by 
$$\forall C\in\B(\RR^2),\ \ \nu_n(C) := 
2\pi\, A_n^{2}\, {\mathbb E}_{(\mu,0)}\big[f(X_n)\, {\bf 1}_{C}(S_n)\big],$$ 
converges weakly to measure $\nu(\cdot) := c_\SS\, (\det\Gamma)^{-1/2}\, \pi(f)\, m_\SS(\cdot)$, 
see \cite{bre}. Since the boundary of $B:=B(s,\varepsilon)$ has zero $\nu-$measure when $\varepsilon\in(0,\varepsilon_{\SS})$, we have: $\lim_{n}\nu_n(B) = \nu(B)$, which is (\ref{LLT-f-bornee-gene}). 
\begin{proof}[Proof of (\ref{h-prrof-pn-lattice}).] Note that 
\begin{enumerate}
\item[-] $\ m_\SS(h) := \int_{\RR^2}h(t)dt$ in Case~{\it (H1)}, 
\item[-] $\ m_\SS(h) := \sum_{n\in\ZZ}\int_{\RR} h(bn,y)\, dy$ in Case~{\it (H2)}, 
\item[-] $\ m_\SS(h) := \sum_{\eta\in \SS}h(\eta)$ in Case~{\it (H3)}. 
\end{enumerate}
Let $h\in{\mathcal H}_2$. Let $\beta$ be a positive real number such that 
Supp($\hat h$)$\, \subset B(0,\beta)$. Without loss of generality, one can suppose that the positive real numbers $\beta$ and $\alpha$ 
(of (\ref{dec-KL})) are such that $\alpha<a/2<\beta$ in Case~{\it (H2)}, 
and $B(0,\alpha)\subset A([-\frac{1}{2},\frac{1}{2}]^2) \subset B(0,\beta)$  
in Case~{\it (H3)}. We set 
\begin{equation} \label{def-K}
K := \big(\overline{B}(0,\beta)\setminus B(0,\alpha)\big)\cap {\mathcal D}.
\end{equation}
Observe that $K$ is a compact subset of $\RR^2\setminus \SS^*$. Let 
$\rho\in(0;1)$ be defined in (\ref{non-ari-spectral-lattice}) w.r.t.~$K$, and set 
$r:=\max(\kappa,\rho)$, where $\kappa$ is defined in (\ref{dec-KL}). 
Using (\ref{dec-KL}) and (\ref{non-ari-spectral-lattice}), we abuse  the notation $O(r^n)$ 
for $Q(u)^n -  \lambda(u)^n\Pi(u)$ when $u\in B(0,\alpha)$, and for $Q(u)^n$ when  
$u\in K$. 
So we have:  
\begin{equation} \label{dec-alpha-beta-lattice}
\forall u\in B(0,\beta),\ \ \ Q(u)^n  =  {\mathbf 1}_{B(0,\alpha)}(u)\, \lambda(u)^n\, \Pi(u) +  O(r^n)
\ \ \mbox{in}\ {\mathcal L}({\mathcal B}),
\end{equation}
where $\Pi(\cdot)$ is the ${\mathcal L}({\mathcal B})$-valued bounded function in (\ref{dec-KL}). Recall that, by hypothesis,  $f\in {\mathcal B}$, $\mu\in(\widehat{\mathcal B})'$. 
Since $\hat h$ is integrable, we then deduce from
Lemma~\ref{fourier-inverse} and (\ref{dec-alpha-beta-lattice}) that
\begin{eqnarray*}
(2\pi)^{2}\, \EE_{(\mu,0)}\big[f(X_n)\, h(S_n)\big] &=& \int_{B(0,\alpha)} P_h(u)\, \lambda(u)^n \mu\big(\Pi(u)f) \, du +  O(r^n) \\
&=& \frac{1}{A_n^2}\, \int_{B(0,\alpha A_n)} P_h\left(\frac{u}{A_n}\right)\, 
 \lambda\left(\frac{u}{A_n}\right)^n \mu\bigg(\Pi\left(\frac{u}{A_n}\right)f\bigg) 
 \, du +  O(r^n). 
\end{eqnarray*}
Next, from (\ref{Dl-2}), $A_n^2\sim nL(A_n)$ and from the fact that $L$ is slowly varying, it can be easily seen that $\lim_n\lambda\big(\frac{u}{A_n}\big)^n = e^{-\langle u,\Gamma u\rangle/2}$. 
Moreover we know by (\ref{b-stat-1}) that 
$\lim_n\, \mu(\Pi(u/A_n)f)  = \pi(f)$. By using (\ref{lambda-ine-CD}), Lebesgue's theorem gives:  
\begin{eqnarray*}
\lim_{n\rightarrow +\infty} \int_{B(0,\alpha A_n)} P_h\left(\frac{u}{A_n}\right)\, 
 \lambda\left(\frac{u}{A_n}\right)^n \mu\bigg(\Pi\left(\frac{u}{A_n}\right)f\bigg) \, du &=&  \pi(f)\, P_h(0)\, \int_{\RR^2} e^{-\langle u,\Gamma u\rangle/2}\, du\\
 &=&  2\pi (\det\Gamma)^{-1/2}\, \pi(f)\, P_h(0).
\end{eqnarray*}
Finally, the Poisson summation formula yields $P_h(0):=\sum_{g\in \SS^*}\hat h(g) = c_\SS\, m_\SS(h)$.
\end{proof}
\subsection{Proof of Proposition~\ref{sec:loi_stable}} \label{proof-sec:loi_stable} 
Proposition~\ref{sec:loi_stable-ter} applied to $f={\bf  1}_{\X}$ gives (\ref{equi-1}). To prove (\ref{equi-2}), let us first state a lemma concerning the sequence $(S_n,S_{n+m}-S_n)_{n,m}$. 
\begin{lem1} \label{couple}
The sequence $(\nu_{n,m})_{n,m}$ of positive measures on $\RR^4$ defined  by 
$$\forall C\in\B(\RR^4),\ \ \nu_{n,m}(C) := (2\pi)^2A_n^2A_m^2\, {\mathbb E}_{(\mu,0)}\big[{\bf 1}_C(S_n,S_{n+m}-S_n)\big],$$
converges weakly, as $\min(n,m)\rightarrow+\infty$, to the measure $\nu$ defined by: $\nu(C) := c_\SS^2\, (det(\Gamma))^{-1}\, m_\SS\otimes m_\SS(C)$.
\end{lem1}
\noindent Before proving this lemma, let us first show how it is used to give
(\ref{equi-2}). Let $T$ be the linear (invertible) endomorphism on $\RR^4$ defined  by: 
$Tw := (u,u+v)$, where we write $w=(u,v)\in\RR^4$, with $u$ and $v$ in $\RR^2$. 
From Lemma~\ref{couple}, the family of measures $\tilde\nu_{n,m}$ 
on $\RR^4$ defined  by 
$$\tilde\nu_{n,m}(C) := \nu_{n,m}({\bf 1}_C\circ T) = (2\pi)^2A_n^2A_m^2\, 
{\mathbb P}_{(\mu,0)}\big((S_n,S_{n+m})\in C\big)$$
converges weakly to $\tilde\nu(C) := \nu({\bf 1}_C\circ T)$ when 
$\min(n,m)\rightarrow+\infty$. But, from Fubini's theorem and since $m_\SS$ is the 
Haar measure, we have $\tilde\nu = \nu$. Now set $B := B(s,\varepsilon)$ for 
$(s,\varepsilon)\in \SS\times(0;\varepsilon_\SS)$. 
Since the boundary of $B\times B$ has zero $\nu-$measure, we obtain the following convergence when $\min(n,m)\rightarrow+\infty$: 
$\lim\tilde\nu_{n,m}(B\times B) = \nu(B\times B) = c_\SS^2\, (det(\Gamma))^{-1}\,m_\SS(B)^2$, which is (\ref{equi-2}). 
\begin{proof}[Proof of Lemma~\ref{couple}] Observe that there exists a continuous $m_{\SS}$-integrable function $h>0$ on 
$\RR^2$ having a compactly supported Fourier transform, see \cite[Section~10.2]{bre}. Define the following function on 
$\RR^4$: $G(w) := h(u)h(v)$, where $w=(u,v)\in\RR^4$, with $u$ and $v$ in $\RR^2$. Then the Fourier 
transform of $G$ is compactly supported on $\RR^4$, and we have $G(w)\, e^{i\langle w,c \rangle}
 = h(u)\, e^{i\langle u,a \rangle}\, h(v)\, e^{i\langle v,b \rangle}$ for any $c=(a,b)\in\RR^4$, 
with $a$ and $b$ in $\RR^2$. 
Therefore, using again classical properties on convergence of positive 
measures \cite{bre},  Lemma~\ref{couple} 
will be established provided that we prove the following: $\forall(h_1,h_2)\in{\mathcal H}_2\times{\mathcal H}_2$, 
\begin{equation}\label{equi-2b}
\lim_{n\rightarrow +\infty} (2\pi)^2A_n^2A_m^2{\mathbb E}_{(\mu,0)}\big[h_1(S_n)\, h_2(S_{n+m}-S_n)\big] 
    =  c_\SS^2\, (det(\Gamma))^{-1}\, m_\SS(h_1)\, m_\SS(h_2). 
\end{equation}
Let $h_1,h_2\in {\mathcal H}_2$ be fixed, and let $\beta>0$ be such that both Supp($\hat h_1$) 
and Supp($\hat h_2$) are contained in $B(0,\beta)$. Real numbers $\alpha,\kappa$ in (\ref{dec-KL}), and $\rho$ in (\ref{non-ari-spectral-lattice}), are chosen 
as in the previous proof, and again we set $r:=\max(\kappa,\rho)$.  
We obtain by using (\ref{dec-alpha-beta-lattice}): $\forall(u,v)\in B(0,\beta)^2$, 
$$Q(u)^n\, Q(v)^m  =  {\mathbf 1}_{B(0,\alpha)}(u)\, {\mathbf 1}_{B(0,\alpha)}(v)\, \lambda(u)^n\, \lambda(v)^m\, 
\Pi(u)\, \Pi(v) +  O(r^{\min(n,m)})\ \ \ \mbox{in}\ {\mathcal L}({\mathcal B}).$$
Using the second formula of Lemma~\ref{fourier-inverse} (with $f={\bf 1}_{\X}$) and  Property~(\ref{b-stat}), the arguments used to prove (\ref{h-prrof-pn-lattice}) can be easily extended to prove (\ref{equi-2b}). 
\end{proof}

%
%
\section{Proof (and complements) for applications \ref{aplli-rho-mix} to \ref{appli-non-stand} of Section~\ref{presentation}} \label{sect-ex-1}

\subsection{$\rho$-mixing case (Proof of Application~\ref{aplli-rho-mix} of Section~\ref{presentation})} \label{subsec-mix}
%
For $p\in \NN^*$ and $q\in \NN^*\cup\{\infty\}$ with $p\le q$, let
${\mathcal G}^{q}_{p}$ denote the $\sigma$-algebra $\sigma(X_p,\ldots,X_q)$ 
generated by $X_p,\ldots,X_q$. The $\rho$-mixing coefficient of $(X_n)_{n\in\NN}$ at horizon $k\geq1$ is defined by  
\begin{equation} \label{Def_rhot}
	\rho(k) := \sup_{j\in\NN^*}\sup\left\{ |\mathrm{Corr}(f;h)|,\ f\in{\mathbb L}^2({\mathcal G}^{j}_{1}), \,  
	h\in {\mathbb L}^2({\mathcal G}_{k+j}^{\infty}) \right\}.
\end{equation}
where $\mathrm{Corr}(f;h)$ is the correlation coefficient of the two random variables $f$ and $g$. 

In Application~\ref{aplli-rho-mix} of Section~\ref{presentation}, the driving Markov chain $(X_n)_{n\in\NN}$ of the MRW $(X_n,S_n)_{n\in\NN}$ is assumed to be $\rho$-mixing, namely  
\begin{equation*}
	\lim_{k\rightarrow +\infty} \rho(k) = 0.
\end{equation*}
The previous property 
is equivalent to the following spectral gap property of the transition kernel $Q$ of $(X_n)_n$ with respect to the Lebesgue space ${\mathbb L}^2(\pi)$, see \cite{rosen}: 
$$\lim_{n\rightarrow +\infty} \sup\big\{\Vert Q^nf-\pi(f)\Vert_{{\mathbb L}^2(\pi)},\ f\in{\mathbb L}^2(\pi),\, \|f\|_2\leq1\big\}=0.$$ 
Classical Markov models satisfying this property are  reviewed in \cite{djl}. 
\begin{proof}[Proof of Application~\ref{aplli-rho-mix} of Section~\ref{presentation}] 
The operator-type hypotheses~{\bf (A1)}-{\bf (A4)} hold with 
${\mathcal B} := {\mathbb L}^2(\pi)$ and $\widehat{\mathcal B} := {\mathbb L}^1(\pi)$: this is established in \cite[Sec.~4-5]{loicsoaz} for additive functionals.  
Extension to general MRW is straightforward. Moreover, since by hypothesis ${\mathbb E}_{(\pi,0)}[|S_1|^2] < \infty$ and ${\mathbb E}_{(\pi,0)}[S_1]=0$, $(S_n/\sqrt n)_n$ converges in distribution under ${\mathbb P}_{(\pi,0)}$ to a Gaussian distribution ${\mathcal N}(0,\Gamma)$, see \cite[Th.~1]{djl}.
Then Application~\ref{aplli-rho-mix} of Section~\ref{presentation} follows from Theorem~III. 
\end{proof} 
Let us mention that the convergence to stable distributions of additive 
functionals associated with $\rho$-mixing Markov chains is investigated in \cite{jara}.
 Unfortunately the non-standard CLT is not studied in \cite{jara}. 

%
%
\subsection{$V$-geometrical ergodicity case (Proof of Application~\ref{appli-v-geo} of Section~\ref{presentation})} \label{sub-v-geo}
%
Given some unbounded function $V : \X\rightarrow  [1,+\infty)$, $(X_n)_{n\in\NN}$ is assumed to be $V$-geometrically ergodic, namely we have $\pi(V) < \infty$ and there exists $\kappa\in(0,1)$ such that we have:
$$\sup_{|f|\leq V}\sup_{x\in \X}\frac{\big|{\mathbb E}_x[f(X_n)] - \pi(f)\big|}{V(x)} = O(\kappa^n),$$
where functions $f : \X\rightarrow\CC$ are assumed to be measurable. The $V$-geometrical ergodicity condition can be investigated with the help of 
the so-called drift conditions. For this fact and for classical examples of such models, 
we refer to \cite{mey}. 
\begin{coro1} \label{result-geo}
Let $\xi : \X\rightarrow \RR^2$ be a $\pi$-centered function taking values in a two-dimensional closed subgroup $\SS$ of $\RR^2$. Set  $S_n:=\sum_{k=1}^n\xi(X_k)$. 

Under Hypotheses~{\bf (A5)}-{\bf (A6)}, if $\vert \xi\vert^\alpha/V$ is bounded for some
$\alpha>0$, then for every initial distribution $\mu$ of $(X_n)_{n\in\NN}$, 
the recurrence set ${\mathcal R}_{(\mu,0)}$ of $(S_n)_n$ satisfies ${\mathcal R}_{(\mu,0)}=\SS$.  
\end{coro1}
\begin{proof}[Proof of Corollary~\ref{result-geo}]
For $\gamma\in(0,1]$, denote by $(\mathcal B_\gamma,\|\cdot\|_\gamma)$ the space of 
measurable $\CC$-valued functions $f$ on $\X$ such that $\|f\|_\gamma =  \sup_{x\in E}\, 
|f(x)|/V(x)^{\gamma}<\infty$. 
Let $\gamma_0\in(0,1)$. From \cite[Sect.~10]{loicsoaz},  
operator-type assumptions {\bf (A1)} to {\bf (A4)}
are fulfilled with  $\mathcal B = \mathcal B_{\gamma_0}$ and 
$\widehat{\mathcal B} = {\mathcal B}_1$ (this is proved in \cite[Lem.~10.1]{loicsoaz} 
with $\widehat{\mathcal B} := {\mathbb L}^1(\pi)$; the case $\widehat{\mathcal B} = 
{\mathcal B}_1$ is similar, use \cite[Lem.~10.4]{loicsoaz}). 
The assumption  of 
Theorem~IV concerning the $\delta_x$'s is obviously fulfilled. Then  
Corollary~\ref{result-geo} follows from 
Theorem~IV.  
\end{proof} 
\begin{proof}[Proof of Application~\ref{appli-v-geo}(a) of Section~\ref{presentation}] 
If $|\xi|/\sqrt V$ is bounded on $\X$, then  the domination assumption of Corollary~\ref{result-geo} is fulfilled, and $(S_n)_n$ satisfies the standard CLT. If moreover $\xi$ is non-sublattice in $\SS$, then the covariance matrix $\Gamma$ of the CLT is automatically positive definite from Remark~\ref{rk-non-matrix}. 
The last remarks together with Corollary~\ref{result-geo} give Application~\ref{appli-v-geo}(a) of Section~\ref{presentation}. 
\end{proof} 

Alternative conditions for the CLT can be found in \cite{jones}. To the best of our knowledge, 
the non-standard CLT has not been investigated for 
$V$-geometrically ergodic Markov chains. 
\begin{proof}[Proof of Application~\ref{appli-v-geo}(b) of Section~\ref{presentation}] 
If $\sup_{x,y} \vert\xi(x,y)\vert^{2+\varepsilon}/(V(x)+V(y)) < \infty$, then {\bf (A1)}-{\bf (A4)} are fulfilled with $\mathcal B = \mathcal B_{\gamma_0}$ (for some $\gamma_0\in(0,1)$) and 
$\widehat{\mathcal B} = {\mathcal B}_1$ (use \cite[Lem.~10.1]{loicsoaz} and \cite[Lemma B.2]{jvl}). 
Moreover, from \cite[Lemma 1]{jvl} and Levy's theorem, it can be easily seen that $S_n:=\sum_{k=1}^n\xi(X_{k-1},X_{k})$ satisfies the standard CLT. We conclude thanks to Theorem~IV.  
\end{proof} 
%
%
%
\subsection{Case of Lipschitz iterative models (Proof of Application~\ref{appli-IFS} of Section~\ref{presentation})} \label{subsec-ifs} 
%
Here $(\X,d)$ is a non-compact metric space in which every 
closed ball is compact.
$\X$ is endowed with its Borel $\sigma$-field ${\mathcal X}$. Let $(V,{\mathcal V})$ be a measurable space, 
let $(\vartheta_n)_{n\geq 1}$ be a  i.i.d.~sequence of random variables taking values in 
$V$, let $F : \X\times V\rightarrow \X$ be a measurable function. Given a $\X$-valued r.v.~$X_0$ independent of $(\vartheta_n)_{n\geq1}$, the random iterative  model associated to 
$((\vartheta_n)_{n\geq 1}, F,X_0)$ is defined by  (see \cite{duf}) 
$$X_n = F(X_{n-1},\vartheta_n),\ \ n\geq 1.$$
Let us consider the two following random variables which are classical in these models (see \cite{duf})~: 
$${\mathcal C} : = \sup\bigg\{\frac{d\big(F(x,\vartheta_1),F(y,\vartheta_1)\big)}{d(x,y)},\ x,y\in \X,\ x\neq y\bigg\}\ \ \ \mbox{and}\ \ \ \ {\mathcal M}  = 1 + {\mathcal C} + d\big(F(x_0,\vartheta_1),x_0\big)$$
where $x_0$ is some fixed point in $\X$. 
\begin{coro1} \label{result-ite}
Assume that ${\mathcal C}<1$ almost surely, that $\EE[{\mathcal M}^{\delta}] < \infty$ for some 
arbitrary small $\delta>0$, and finally that $\xi$ satisfies the following weighted-Lipschitz 
condition: 
\begin{equation} \label{xi-iterative} 
\exists C\geq 0,\ \exists b\geq0,\ \forall (x,y)\in \X^2,\ \ \big|\xi(x)-\xi(y)\big| \leq C\, d(x,y)\, \big[1+d(x,x_0)+d(y,x_0)\big]^b,
\end{equation} 
Then, under Hypotheses~{\bf (A5)}-{\bf (A6)}, the AF $(S_n)_n$ defined in (\ref{add-fct}) satisfies 
${\mathcal R}_{(\mu,0)}=\SS$ for every initial distribution $\mu$ on $\X$. 
\end{coro1} 
\begin{proof} 
Set: $\forall x\in \X,\ p(x) = 1+ \, d(x,x_0)$. For $0<\alpha\leq 1$ and $\gamma>0$, let ${\mathcal B}_{\alpha,\gamma}$ be the space of all $\CC$-valued functions on $\X$ satisfying the following condition
$$m_{\alpha,\gamma}(f) =  
\sup\bigg\{\frac{|f(x)-f(y)|}{d(x,y)^\alpha\, p(x)^{\alpha\gamma}\, p(y)^{\alpha\gamma}},\ x,y\in \X,\  
x\neq y\bigg\}\, <\, +\infty.$$
Set  $|f|_{\alpha,\gamma}  = \sup_{x\in \X}\, |f(x)|/p(x)^{\alpha(\gamma+1)}$ and
$\|f\|_{\alpha,\gamma} = m_{\alpha,\gamma}(f) + |f|_{\alpha,\gamma}$. Then
$({\mathcal B}_{\alpha,\gamma},\|\cdot\|_{\alpha,\gamma})$  is a Banach space 
(this corresponds to the space ${\mathcal B}_{\alpha,\beta,\gamma}$ of 
\cite[Sect.~11]{loicsoaz} in case $\beta=\gamma$). 
Now let us assume that $\gamma>b+1$ and $2\alpha\gamma < \delta$, and consider $\gamma'>\gamma$ such that $b+1+(\gamma'-\gamma)\leq \gamma$ and $\alpha(\gamma'+\gamma)\leq\delta$. 
Using assumptions on ${\mathcal C}$ and ${\mathcal M}$, it follows from \cite[p.~483]{loicsoaz} that the operator-type assumptions 
{\bf (A1)} to {\bf (A4)} are fulfilled with $\mathcal B := {\mathcal B}_{\alpha,\gamma}$ and with $\widehat{\mathcal B}$ defined as 
the Banach space of all the $\CC$-valued functions on $\X$ satisfying $|f|_{\alpha,\gamma'}  = \sup_{x\in \X}\, |f(x)|/p(x)^{\alpha(\gamma'+1)}$. 
Note that inclusion $\widehat{\mathcal B} \subset {\mathbb L}^1(\pi)$ is continuous since $p(\cdot)^{\alpha(\gamma'+1)}$ is $\pi$-integrable (use $\alpha(\gamma'+1)\leq\delta$ and \cite[Prop.~11.1]{loicsoaz}). Then  Corollary~\ref{result-ite} follows from Theorem~IV.  
\end{proof}
\begin{proof}[Proof of Application~\ref{appli-IFS} of Section~\ref{presentation}] 
From ${\mathcal C}<1$ and $\EE[\, {\mathcal M}^{2(b+1)}] < \infty$, $(S_n)_n$ satisfies the standard CLT, see \cite[Prop.~11.3]{loicsoaz}. From non-sublattice condition in $\SS$, the covariance matrix $\Gamma$ of the CLT is positive definite, see Remark~\ref{rk-non-matrix}. We conclude thanks to Corollary~\ref{result-ite}. 
\end{proof}

A more precise use of results of  \cite[Sect.~11]{loicsoaz} enables to obtain 
the conclusion of Corollary~\ref{result-ite} under some weaker mean contractive conditions on ${\mathcal C}$ 
(instead of ${\mathcal C}<1$ a.s.). Alternative conditions for the CLT can be found in \cite{duf,benda,wu-shao} and references therein. 
%
%
\subsection{Complement on Application~\ref{appli-non-stand} of Section~\ref{presentation}}  \label{subsec-aff-rec}
Let $(B_n,A_n)_n$ be a sequence of i.i.d.~random variables with values
in ${\mathbb R}^2\rtimes Sim({\mathbb R}^2)$, independent of $X_0$, 
where $Sim({\mathbb R}^2)$ is the 
similarity group of ${\mathbb R}^2$. Let us consider the affine 
iterative model
$$X_n=A_nX_{n-1}+B_n.$$
Now assumptions and statements of Application~\ref{appli-non-stand} (Section~\ref{presentation}) are specified.  
Let $\nu$ denote the distribution of $(B_1,A_1)$. 
We write $\bar\nu$ the projection of $\nu$ on $Sim({\mathbb R}^2)$, we denote
by $G_{\bar\nu}$ the closed subgroup of $Sim({\mathbb R}^2)$ generated by the
support of $\bar\nu$. We recall that $Sim({\mathbb R}^2)={\mathbb R}^*_+\times
 O({\mathbb R}^2)$, where $ O({\mathbb R}^2)$ is the orthogonal group. We also denote by 
$|\cdot|$ the matrix norm associated with the euclidean norm on $\RR^2$. 
We suppose that:
\begin{itemize}
\item there exists a unique stationary distribution $\pi$ and its support is
unbounded, 
\item no affine subspace of ${\mathbb R}^2$ is invariant by the support of $\nu$, 
\item ${\mathbb E}[|A_1|^2]= 1$, ${\mathbb E}[\vert A_1|^2\log|A_1|]<\infty$
and ${\mathbb E}[|B_1|^2]<\infty$, 
\item the projection of $G_{\bar\nu}$ on ${\mathbb R}^*_+$ is equal
to ${\mathbb R}^*_+$.
\end{itemize}
In this case, measure $\pi$ admits an expectation in $\RR^2$, called $m_0$, 
and there exists a gaussian random variable $Z$ such that, for every $x\in{\mathbb R}^2$, 
under the probability measure ${\mathbb P}_x$, the following sequence of random variables 
$$\bigg(\frac 1{\sqrt{n\log(n)}}\big(\sum_{k=1}^nX_k-nm_0\big)\bigg)_{n\ge 1}$$
converges in distribution to $Z$ (\cite[Th. 1.5]{BDG}). 
\begin{coro1} \label{result-affine} 
If $x\mapsto x$ is non-sublattice in ${\mathbb R}^2$, then, for every $(x,s)\in{\mathbb R}^2\times{\mathbb R}^2$, 
we have 
$${\mathbb P}_x\bigg(\forall \varepsilon>0,\ \bigg\vert\sum_{k=1}^nX_k-nm_0-s\bigg\vert
        <\varepsilon,\ \text{i.o.}\bigg)=1.$$
\end{coro1}
\begin{proof}
Again we apply Theorem~IV. Thanks to Lemmas 3.9 and 3.12 of \cite{BDG}, operator-type assumptions {\bf (A1)} to {\bf (A4)} hold true on some Lipschitz weighted spaces similar to those introduced in Subsection~\ref{subsec-ifs}.  
Probabilistic-type assumptions follow from the above non-standard CLT 
(see also \cite[Prop 3.18]{BDG}) and from the non-sublattice assumption on $\xi(x) = x$.  
\end{proof}


\begin{appendix}
\section{Complement on Assertion~$(a)$ in Proposition~\ref{sec:KL}} \label{conv-gauss-mambda} 
\begin{prop1} \label{lem-DL2}
Assume that $Q$ is strongly ergodic on $\mathcal B$ (see {\bf (A1)}) and that 
\begin{enumerate}
  \item[{\bf (A2')}]
There exists $\alpha>0$ such that, for every $t\in B(0,\alpha)$, we have $Q(t)\in\mathcal L(\mathcal B)\cap{\mathcal L}(\widehat{\mathcal B})$ and: 
\begin{equation} \label{cont-kel-liv} 
\forall t\in B(0,\alpha),\ \lim_{h\rightarrow 0}\|Q(t+h)-Q(t)\|_{{\mathcal B},\widehat{\mathcal B}}=0,
\end{equation}
  \item[{\bf (A3')}]
There exist $\kappa_1\in[0,1)$ and $C\in(0;+\infty)$ such that 
$$\forall n\geq1,\ \forall t\in B(0,\alpha),\ \forall f\in\mathcal B, \ \
\|Q(t)^nf\|_{\mathcal B} \leq C\, \kappa_1^n\, \|f\|_{\mathcal B} + C\, \|f\|_{\widehat{\mathcal B}}.$$  
\end{enumerate}
Then Property~(\ref{dec-KL}) of Proposition~\ref{sec:KL} is fulfilled. Moreover properties (\ref{Dl-2}) and {{\bf (A5)}} are equivalent (with the same covariance matrix $\Gamma$ and the 
same function $L(\cdot)$).  
\end{prop1} 
\noindent 
Since {{\bf (A2)}}-{{\bf (A3)}} imply {\bf (A2')}-{\bf (A3')}, Proposition~\ref{lem-DL2} completes the proof of Proposition~\ref{sec:KL}$(a)$. 
\begin{proof} 
The fact that (\ref{dec-KL}) holds under Hypotheses~ {\bf (A1)} and {\bf (A2')}-{\bf (A3')} follows from \cite[p.~428]{ihp1}. The equivalence between (\ref{Dl-2}) and {{\bf (A5)}} is proved in \cite[lem.~4.2]{ihp1} when $A_n=\sqrt n$ in {{\bf (A5)}}, see also \cite[Lem.~5.2]{loicsoaz}. 
When the non-standard CLT holds in {{\bf (A5)}}, the proof is similar, we just outline below the main arguments. 
Without loss of generality, we suppose that $\Gamma$ is the identity matrix. First observe that we have by Lemma~\ref{fourier} (applied with $f={\bf 1}_\X$ and $m=0$):
\begin{equation} \label{charac-clt}
\forall u\in B(0,\alpha),\ \ \EE_{(\pi,0)}\big[e^{i\langle u,S_n\rangle}\big] = \pi\big(Q(u)^n {\bf 1}_\X\big).
\end{equation}
The proof of the "if-part" in Proposition~\ref{lem-DL2} is easy: indeed, assume that (\ref{Dl-2}) holds, and let $(A_n)_n$ be a sequence of positive real numbers such that  $A_n^2\sim nL(A_n)$. From (\ref{charac-clt}) and (\ref{dec-KL}), we obtain for any fixed $t\in\RR^2$ and for $n$ sufficiently large 
$$\EE_{(\pi,0)}\big[e^{i\langle t,S_n/A_n\rangle}\big] = 
\lambda\big(t/A_n\big)^n\, \pi\big(\Pi(t/A_n){\bf 1}_\X\big) + O(\kappa^n).$$ 
Using (\ref{Dl-2}), $A_n^2\sim nL(A_n)$ and the fact that $L$ is slowly varying, one can easily see that $\lim_n \lambda\big(t/A_n\big)^n = e^{-|t|^2/2}$. Hence the desired CLT in Hypothesis~{{\bf (A5)}} holds true.  

\noindent Conversely, assume that {{\bf (A5)}} holds. Let us prove that  the function $\lambda(\cdot)$ in (\ref{dec-KL}) satisfies: 
$$\psi(u) := \frac{\lambda(u)-1}{|u|^2\, L\big(|u|^{-1}\big)} + 1/2 \rightarrow 0\ \ \ \mbox{when}\ u\rightarrow 0.$$
{}From Levy's theorem, we have: $\forall t\in B(0,\alpha),\ \lim_n\EE_{(\pi,0)}
[e^{i\langle t, S_n/A_n\rangle}] = \exp(-|t|^2/2)$. Thus, by using 
(\ref{charac-clt}), (\ref{dec-KL}) and the complex logarithm function $\log(\cdot)$, this gives 
$\lim_n\, n\log\lambda(t/A_n) = -|t|^2/2$, from which we easily deduce: $\lim_n \psi(t/A_n)=0$ 
(use $\log(z) \sim (z-1)$ when $z\rightarrow 1$, $n\sim A_n^2/L(A_n)$ and $L(A_n) \sim L(A_n/|t|)$). 
More precisely, by using the classical refinement of Levy's theorem in terms of uniform convergence on compact sets and 
the fact that the property 
$\lim_{x\rightarrow +\infty}\frac {L(kx)} {L(x)}=1$ is uniform in $k$ 
on each compact subset of $(0;+\infty)$ (according to formula (\ref{karamata})), 
one can see that the limit $\lim_n \psi(t/A_n)=0$ is uniform on $C := \{t\in\RR^2 : \alpha/2\leq |t| < \alpha\}$ (see \cite[lem.~4.2]{ihp1} for details). So,  given $\varepsilon>0$, one can choose $N_0=N_0(\varepsilon)$ such that: $n\geq N_0,\, t\in C\ \Rightarrow\ |\psi(t/A_n)|<\varepsilon$. 
Next, since $\lim_n A_{n+1}/A_n=1$, one can suppose that $N_0$ is such that: 
$\forall n\geq N_0,\ 1/(2A_n) < 1/A_{n+1}$. From that, we easily deduce that 
$\cup_{n\geq N_0} C/A_n = \{t\in\RR^2 : 0 < |u| < \alpha/A_{N_0}\}$. Therefore 
we have: $0 < |u| < \alpha/A_{N_0}\ \Rightarrow\ |\psi(u)|<\varepsilon$. 
\end{proof}
\begin{rqe1}
Expansion~(\ref{Dl-2}) may be adapted to cover the convergence in distribution of $S_n$ (properly normalized) to a stable 
distribution of index $0<p<2$. Then Propositions~\ref{sec:loi_stable-ter} and \ref{sec:loi_stable} extend  
(with $A_n$ such that $A_n^p \sim n\, L(A_n)$), but we have $\sum_{n\geq1} A_n^{-2}<\infty$
since for $n$ large enough,  $A_n^{p/2-1}\le L(A_n)$,
$nL(A_n)\le 2A_n^p$ 
and thus $A_n^{-2}\le  (2/n)^{4/(p+2)}$. 
Therefore we obtain: $\sum_{n\geq1}{\mathbb P}_{(\mu,0)}\big( \vert S_n-s\vert<\varepsilon\big) < \infty$. This gives the  expected transience property. 
\end{rqe1}
\section{Complement on Assertion~$(c)$ in Proposition~\ref{sec:KL}} \label{app-B}

Without loss of generality, we suppose that the MRW $(X_n,S_n)_{n\in{\mathbb N}}$ 
is the canonical version defined on $\Omega = (\X\times{\mathbb R}^2)^{\mathbb N}$.
We recall that, in the sense given in (\ref{Sn-supported-by-H}), $S_1$ takes its values
in a closed two-dimensional subgroup $\SS$ of ${\mathbb R}^2$: this corresponds to 
cases {\it (H1) (H2) (H3)} described at the beginning of Section~\ref{simp-rks-S}. 
In this appendix, we prove that, under Hypotheses~{{\bf (A1)}}-{{\bf (A4)}}, 
Property~(\ref{non-ari-spectral-lattice}) 
is linked to the non-sublattice condition of Hypothesis~{\bf (A6)}. 
To that effect, introduce the following: 
\begin{defi1} \label{def-ari-sublat}
Under Hypothesis~(\ref{Sn-supported-by-H}) we shall say that 
$(S_n)_n$ is {\bf arithmetic} in $\SS$ w.r.t.~$\mathcal B$ if there exist 
$t\in{\mathbb R}^2\setminus \SS^*$,
$\lambda\in\mathbb C$ with $|\lambda|=1$, $w\in\mathcal B$ such that, for $\pi$-almost every $x\in\X$, 
we have $|w(x)|=1$ and the following property: 
\begin{equation} \label{fle-non-ari}
\forall n\geq 1,\ \ e^{i \langle t , S_n \rangle}\, w(X_n) = \lambda^n\, w(x)\ \ {\mathbb P}_{(x,0)}-\mbox{a.s.}. 
\end{equation}
\end{defi1}

Under Hypothesis~{{\bf (A1)}}-{{\bf (A4)}}, we consider the set 
$$G:=\{t\in{\mathbb R}^2\ :\ r(Q(t))=1\}.$$
Recall that the dual subgroup $\SS^*$ of $\SS$ is defined in (\ref{H*}). Since $Q(0) = Q$, $r(Q)=1$ and $Q(\cdot)$ is $\SS^*$-periodic, $\SS^*$ is contained in $G$. 
\begin{prop1} \label{prop-red-lattice}
Assume that Hypotheses~{{\bf (A1)}}-{{\bf (A4)}} hold true. Then the following assertions hold: \\[0.12cm] 
(i) Property~(\ref{non-ari-spectral-lattice})  $\Leftrightarrow\ G=\SS^*\ \Leftrightarrow\ (S_n)_n$ is not arithmetic in $\SS$ w.r.t. $\mathcal B$; \\[0.12cm]
(ii) If Hypotheses~{{\bf (A5)}}-{{\bf (A6)}} hold, then Property~(\ref{non-ari-spectral-lattice})  is fulfilled; \\[0.12cm]
(iii) If $(S_n)_n$ is sublattice in $\SS$ and the function $\chi(\cdot)$ in 
(\ref{fct-beta}) is such that, for every $t\in{\mathbb R}^2$, we have 
$e^{i\langle t,\chi(\cdot)\rangle}\in\mathcal B$ , then Property~(\ref{non-ari-spectral-lattice})  does not hold. 
\end{prop1}

\begin{rqe1} \label{rk-sub-lat-ifs}
In $\rho$-mixing or $V$-geometrical ergodicity cases (see Subections~\ref{subsec-mix} and \ref{sub-v-geo}), 
the condition on $\chi(\cdot)$ in Assertion~{\it (iii)} is automatically fulfilled, so that the non-sublattice assumption is equivalent to Condition~(\ref{non-ari-spectral-lattice}). For Lipschitz iterative models, 
the non-sublattice assumption is just a sufficient condition
for Condition~(\ref{non-ari-spectral-lattice})  to hold true on the weighted-Lipschitz spaces defined in Subsection~\ref{subsec-ifs}  (because the condition on $\chi(\cdot)$ 
in Assertion~{\it (iii)} of Proposition~\ref{prop-red-lattice} is not automatically fulfilled). 
The non-arithmeticity condition, which is equivalent to Condition~(\ref{non-ari-spectral-lattice}), 
can be simplified in the special case of additional 
functionals (see \cite[Section 5]{loicsoaz}).
\end{rqe1}

When $\SS=\RR^2$ and $S_n$ is an additive functional (see (\ref{add-fct})), Proposition~\ref{prop-red-lattice} is established in \cite[Section 12]{loicsoaz}. 
Here we give the adaptation to general MRWs and subgroups $\SS$. 

\begin{proof}[Proof of Proposition~\ref{prop-red-lattice}]  
First, using Hypotheses~{{\bf (A1)}}-{{\bf (A4)}} and the same arguments as in \cite[Lem.~12.1]{loicsoaz}, we obtain 
\begin{equation} \label{rQt-strict-1}
\forall t\in\RR^2\setminus G,\ \  r(Q(t))<1.
\end{equation}
Second, an easy adaptation of \cite[Lem.~12.3]{loicsoaz} shows that, 
for any compact subset $K$ of $\RR^2\setminus G$, there exists $\rho=\rho(K)\in[0,1)$ 
such that 
\begin{equation} \label{non-ari-spectral-G} 
\sup_{t\in K}\Vert Q(t)^n\Vert_{\mathcal B} = O(\rho^n)
\end{equation}
(consider compact subsets $K$ of
${\mathbb R}^2\setminus G$ instead of compact subsets of ${\mathbb R}^2\setminus\{0\}$ in the 
proof of \cite[Lem.~12.3]{loicsoaz}). 
Next, since $\SS^*\subset G$, the previous property yields the first 
equivalence in {\it (i)}: indeed, if $G=\SS^*$, then 
(\ref{non-ari-spectral-G}) obviously gives (\ref{non-ari-spectral-lattice}) . Conversely, if (\ref{non-ari-spectral-lattice})  is true, then 
for every $t\in\RR^2\setminus \SS^*$ we have $r(Q(t))<1$, thus $t\in\RR^2\setminus G$. Therefore  
Condition~(\ref{non-ari-spectral-lattice})  gives $G\subset \SS^*$, hence $G=\SS^*$. 

\noindent In addition \cite[Lem.~12.1]{loicsoaz} gives the following equivalence: 

\noindent{\it  Property~(A): we have $t\in G$ if and only if there exist $\lambda\in{\mathbb C}$, $|\lambda|=1$, 
and $w\in{\mathcal B}$, $w\neq 0$, such that we have the following equality: 
$Q(t)w=\lambda\, w \ \ \mbox{in}\ {\mathcal B}$. Moreover the previous function $w(\cdot)$ is such that $|w| = \pi(|w|)\, $ $\pi$-a.s.. 
}

\noindent The fact that $Q(t)w=\lambda\, w$ implies $|w| = \pi(|w|)\, $ $\pi$-a.s. is easy to obtain. Indeed, we have $|w| = |\lambda^n\, w| = |Q(t)^nw| \leq Q^n|w|$ for every $n\geq1$,  thus we deduce from {\bf (A4)} that $|w| \leq \pi(|w|)\, $ $\pi$-a.s. So $g := \pi(|w|) - |w|$ is nonnegative, and $\pi(g)=0$, hence $|w| = \pi(|w|)\, $ $\pi$-a.s.. 

\noindent Finally, from the previous property, we can deduce the following. 

\noindent {\it Property~(B): we have $\SS^*\neq G$ if and only if there exist $t\in \RR^2\setminus \SS^*$, $\lambda\in{\mathbb C}$, $|\lambda|=1$, and $w\in{\mathcal B}$, $w\neq 0$, such that $|w|=1$ $\pi$-a.s. and $Q(t)w=\lambda\, w$ in ${\mathcal B}$. }

\noindent To prove the second equivalence in {\it (i)}, one needs the following. 
\begin{lem1} \label{fct-propre-vers-E}
We have $\SS^*\neq G$ if and only if there exist $t\in \RR^2\setminus \SS^*$, $\lambda\in{\mathbb C}$, $|\lambda|=1$, 
and $w\in{\mathcal B}$, $w\neq 0$, such that for $\pi$-a.e.~$x\in\X$ we have $|w(x)|=1$ and 
\begin{equation} \label{eq-w}
\forall n\geq 1,\ \ {\mathbb E}_{(x,0)}\big[e^{i \langle t , S_n \rangle}\, w(X_n)\big] = \lambda^n\, w(x). 
\end{equation}
\end{lem1}
\begin{proof}
Assume that $\SS^*\neq G$, and let $(t,\lambda,w)$ be as stated in Property~(B). 
Then we have: $\forall n\geq 1,\ Q(t)^nw=\lambda^n\, w$ in ${\mathcal B}$. 
Since, by hypothesis, ${\mathcal B}\subset{\mathbb L}^1(\pi)$ with continuous inclusion, 
it follows that $Q(t)^nw=\lambda^n\, w$ in ${\mathbb L}^1(\pi)$, 
hence we have (\ref{eq-w}) for $\pi$-a.e.~$x\in\X$ (use Lemma~\ref{fourier} with $m=0$).  
Conversely, let $t\in \RR^2\setminus \SS^*$ and $(\lambda,w)$ as stated in 
Lemma~\ref{fct-propre-vers-E}. Then we have for $\pi$-a.e.~$x\in\X$: 
$\forall n\geq1,\ Q(t)^nw(x)=\lambda^n\, w(x)$. This implies that $t\in G$. 
Indeed, if $t\notin G$, then by (\ref{rQt-strict-1}) we would have  $r(Q(t))<1$, 
thus $\lim_nQ(t)^nw=0$ in ${\mathcal B}$, and so in ${\mathbb L}^1(\pi)$: this 
would give the property: $w=0$ $\pi$-a.s., which is impossible since by hypothesis 
$|w|=1$ $\pi$-a.s..
\end{proof}
\noindent Using the facts that ${\mathbb P}_{(x,0)}$ is a probability measure and $|w| = 1\, $ $\pi$-a.s., the property stated in Lemma~\ref{fct-propre-vers-E} is equivalent to the arithmeticity of $(S_n)_n$ in $\SS$ w.r.t.~$\mathcal B$, which proves 
the second equivalence in {\it (i)}. 

\noindent Now we prove Assertion~{\it (ii)} of Proposition~\ref{prop-red-lattice}. 
Under Hypothesis~{{\bf (A1)}}-{{\bf (A4)}}, $G$ is a closed subgroup of $\RR^2$, 
and under the additional Hypothesis~{{\bf (A5)}}, $G$ is discrete, see~\cite[Prop.~12.4]{loicsoaz}. 
Observe that, since $\SS^*\subset G$, we have $G^*\subset \SS$. 
To prove Assertion~{\it (ii)} of Proposition~\ref{prop-red-lattice}, one needs to use the following statement, which is an easy adaptation of the proof of \cite[Prop.~12.4]{loicsoaz}: 

\noindent{\it Property~(C): 
there exist a bounded measurable function $\chi:{\mathbb X}
\rightarrow {\mathbb R}^2$ and a family $(\beta_t)_{t\in G}$ of real numbers
such that, for $\pi$-almost every $x\in\X$, we have }
\begin{equation} \label{condC}
\forall t\in G,\ \forall n\geq1,\ \big\langle t, S_n + \chi(X_n) - \chi(x)\big\rangle \in 
n \beta_t + 2\pi\Z\ \ {\mathbb P}_{(x,0)}-a.s.. 
\end{equation}
The fact that $G$ is discrete plays an important role in Property~(C) to obtain the 
existence of the above function $\chi$, which does not depend on $t$.   

\noindent 
Assume that Condition~(\ref{non-ari-spectral-lattice})  is not fulfilled. Then from 
Assertion~{\it (i)} of Proposition  \ref{prop-red-lattice}, 
$\SS^*$ is a proper subgroup of $G$. Hence $G^*$ is a proper subgroup of $\SS$. 
Consequently, from Property~(C), $(S_n)_n$ is sublattice in $\SS$. This proves {\it (ii)}. 

\noindent Finally we establish  Assertion~{\it (iii)} of Proposition~\ref{prop-red-lattice}. Suppose that $(S_n)_n$ is  
sublattice in $\SS$, with $\SS_0$, $\chi(\cdot)$ and $(\beta_t)_{t\in \SS_0^*}$ as indicated in (\ref{fct-beta}), and with the additional condition: $\forall t\in\RR^2,\ e^{i\langle t,
\chi(\cdot)\rangle}\in\mathcal B$. Since by hypothesis $\SS_0$ is strictly contained 
in $\SS$, there exists $t_0\in \SS_0^*\setminus \SS^*$. We deduce from (\ref{fct-beta}) 
that, for $\pi$-almost every $x\in\X$, we have 
$$\forall n\geq1,\ \ e^{i\langle t_0,S_n\rangle} e^{i\langle t_0,\chi(X_n)\rangle} = 
   e^{i n \beta_{t_0}} e^{i\langle t_0,\chi(x)\rangle}\ \ {\mathbb P}_{(x,0)}-a.s..$$
So we obtain (\ref{fle-non-ari}) with $\lambda:=e^{i\beta_{t_0}}$ and 
$w:=e^{i\langle t_0,\chi(\cdot)\rangle}\in\mathcal B$. Hence, from Assertion~(i) of Proposition~\ref{prop-red-lattice}, Condition~(\ref{non-ari-spectral-lattice})  does not hold. 
\end{proof}
\begin{rqe1} \label{dern-rque}
When $(S_n)_n$ is an AF (see (\ref{add-fct})), Conditions (\ref{fle-non-ari}) and (\ref{fct-beta}) in Definition~\ref{def-ari-sublat} may be stated only for $n=1$ and specified with absorbing sets (instead of properties fulfilled $\pi$-a.s.), see \cite{loicsoaz} and Remark~\ref{sec:CLT-xi}. Similarly, for general MRW, if $\mathcal B$ is composed of $\pi$-classes of functions (for instance $\mathcal B ={\mathbb L}^2(\pi)$), then the equivalence in Lemma~\ref{fct-propre-vers-E} is valid when (\ref{eq-w}) holds for $n=1$. Indeed this condition says that for $\pi$-a.e.~$x\in\X$: $Q(t)w(x)=\lambda\, w(x)$. So $Q(t)w=\lambda\, w$ in $\mathcal B$ and the proof of Lemma~\ref{fct-propre-vers-E} can be then repeated. Consequently, under the previous condition on $\mathcal B$, Conditions (\ref{fle-non-ari}) and (\ref{fct-beta}) may be  also stated only for $n=1$ (and then (\ref{gamma-t}) is not relevant). \\
\end{rqe1}
\end{appendix}

\end{document}